%% file: hhs-arXiv.tex
\documentclass[smallextended,numbook,envcountsame]{svjour3}
\smartqed
\usepackage{graphicx}
\usepackage[english]{babel}
\usepackage{amsmath}
\usepackage{amsfonts}
\usepackage{color}
\usepackage{enumitem}
\usepackage{tocloft}

\newcommand{\Ad}{\textup{Ad}}
\renewcommand{\Re}{\textup{Re }}
\newcommand{\id}{\textup{id}}

\newcommand {\rr}{{\mathbb{R}}}
\newcommand {\cc}{{\mathbb{C}}}
\newcommand {\zz}{{\mathbb{Z}}}
\newcommand {\nn}{{\mathbb{N}}}
\newcommand {\g}{{\mathfrak{g}}}

\newcommand {\la}{{\mathfrak{a}}}
\newcommand {\lac}{{\mathfrak{a}_{\cc}}}
\newcommand{\laos}{\la^{*(2)}}
\newcommand{\laosplus}{(\la_+^*)^2}
\newcommand{\lambdaisom}{\Lambda}
\newcommand{\projGM}{\tilde\sigma}
\newcommand {\laregs}{{\mathfrak{a}}_{\mathrm{reg}}^*}
\newcommand {\lp}{{\mathfrak{p}}}

\newcommand {\lk}{{\mathfrak{k}}}

\newcommand {\lnn}{{\mathfrak{n}}}

\renewcommand{\phi}{\varphi}
\renewcommand{\epsilon}{\varepsilon}
\newcommand {\1}{{\mathbf{1}}}

\DeclareMathOperator{\Ind}{Ind}
\DeclareMathOperator{\ad}{ad}

\DeclareMathOperator{\Op}{Op}
\DeclareMathOperator{\supp}{supp}
\DeclareMathOperator{\sign}{sign}

\newcommand{\DD}{\mathbb{D}}
\newcommand{\RE}{\operatorname{Re}}
\newcommand{\IM}{\operatorname{Im}}
\newcommand{\signature}{\operatorname{sgn}}
\newcommand{\bigoh}{\mathcal{O}}
\newcommand{\Cinfty}{\ensuremath{C^{\infty}}}
\newcommand{\Ccinfty}{\ensuremath{C_c^{\infty}}}
\newcommand{\Dprime}{\ensuremath{\mathcal{D}^\prime}}
\newcommand{\intd}{\mspace{0.5mu}\operatorname{d}\mspace{-2.5mu}}
\newcommand{\intdbar}{\mspace{0.5mu}\operatorname{d}\mspace{-10mu}{}^-\!}

\newcommand{\tofrom}[2]{[#1\leftarrow #2]} 
\newcommand{\XG}{X_\Gamma}

\newcommand{\psdiff}{pseudo-differential}

\numberwithin{equation}{section}

\begin{document}

\title{Patterson--Sullivan distributions in higher rank}

\author{S.~Hansen \and J.~Hilgert\footnote{Part of this research was done at the Hausdorff Research Institute for Mathematics in the context of the trimester program ``Interaction of Representation Theory with Geometry and Combinatorics''} \and M.~Schr\"{o}der
\footnote{
Partially supported by the DFG-IRTG 1133
``Geometry and Analysis of Symmetries''}
}
\institute{
\textsc{S\"onke Hansen}
\and
\textsc{Joachim Hilgert}
\and
\textsc{Michael Schr\"oder}\\
\texttt{{soenke.hansen@math.upb.de $\cdot$ hilgert@math.upb.de $\cdot$ michaoe@math.upb.de}}
\at
\textsc{Institut f\"ur Mathematik, Universit\"at Paderborn, Warburger Str. 100, 33098 Paderborn, Germany.}
}

\date{Version: May 29, 2011}

\maketitle

\begin{abstract} For a compact locally symmetric space $\XG$ of non-positive curvature, we consider sequences of normalized joint eigenfunctions which belong to the principal spectrum of the algebra of invariant differential operators. Using an $h$-\psdiff\ calculus on $\XG$, we define and study lifted quantum limits as weak$^*$-limit points of Wigner distributions.
The Helgason boundary values of the eigenfunctions allow us to construct Patterson--Sullivan distributions on the space of Weyl chambers. These distributions are asymptotic to lifted quantum limits and satisfy additional invariance properties, which makes them useful in the context of quantum ergodicity. Our results generalize results for compact hyperbolic surfaces obtained by Anantharaman and Zelditch.
\keywords{Patterson--Sullivan distributions \and Wigner distributions \and quantum ergodicity \and lifted quantum limits \and locally symmetric spaces \and geometric pseudo-differential analysis \and Weyl chamber flow} \subclass{53C35 \and 58C40 \and 58J50}
\end{abstract}


\input hhs-intro	
\input hhs-geom		
\input hhs-hebv		
\input hhs-psd		
\input hhs-oint		
\input hhs-qlim		
\input hhs-asymp	

\begin{acknowledgements}
We thank J.~M\"{o}llers, A.~Pasquale, and in particular N.~Anantharaman and S.~Zelditch for helpful discussions.
Special thanks go to M. Olbrich for providing the idea of a proof of Proposition \ref{Proposition Olbrich}.
\end{acknowledgements}

\end{document}

%% file: hhs-intro.tex

\section{Introduction}

For a locally symmetric space $\XG$ of non-positive curvature,
we consider sequences, $(\varphi_h)_h\subset L^2(\XG)$,
of normalized joint eigenfunctions which belong to the principal spectrum
of the algebra of invariant differential operators.
Using a $h$-\psdiff\ calculus on $\XG$, we define and study \emph{lifted quantum limits}
or microlocal lifts as weak$*$-limit points of Wigner distributions
\[ W_h:a\mapsto \big(\Op_{h}(a)\varphi_h\mid\varphi_h\big)_{L^2(\XG)}. \]
Here, $h^{-1}$ is the norm of a spectral parameter associated with $\varphi_h$,
and $h\downarrow 0$ through a strictly decreasing null sequence.
Lifted quantum limits are positive Radon measures supported in the cosphere bundle.
The problem of quantum ergodicity asks for a description of the lifted quantum limits.
Using the Helgason boundary values of the $\varphi_h$,
we construct \emph{Patterson--Sullivan distributions} on the space 
of Weyl chambers.
In the context of quantum ergodicity, Patterson--Sullivan distributions are important
because they are asymptotic to lifted quantum limits and satisfy invariance properties.

For compact hyperbolic surfaces $X_{\Gamma}=\Gamma\backslash\mathbb{H}$,
the asymptotic equivalence of lifted Wigner distributions and Patterson--Sullivan distributions
was observed by Anantharaman and Zelditch \cite{AZ}.
While it was known from earlier work (see \cite{Z89,Wol}) that lifted quantum limits
on compact hyperbolic surfaces are invariant under geodesic flows it turned out that
Patterson--Sullivan distributions are themselves invariant under the geodesic flow.
Moreover, in \cite{AZ} it is shown that they have an interpretation in terms of
dynamical zeta functions which can be defined completely in terms of the geodesic flow.

Although lifted quantum limits do not depend on the specific \psdiff\ calculus
chosen for their definition, it is useful, for establishing invariance properties,
to have an equivariant calculus.
For hyperbolic surfaces,
based on the non-euclidean Fourier analysis 
and closely following the euclidean model,
such a calculus was provided by  Zelditch \cite{Z84}.
In \cite{S} this calculus was extended to rank one symmetric spaces.
Using this calculus the construction of the Patterson--Sullivan distributions and
the proof of the asymptotic equivalence from \cite{AZ} has been generalized in \cite{HS09}.
However, due to singularities arising from Weyl group invariance,
it is difficult to construct an equivariant non-euclidean \psdiff\ calculus in higher rank; see \cite{S}.
Silberman and Venkatesh \cite{SV07}, generalizing work of Zelditch and Wolpert for surfaces
to compact locally symmetric spaces,
introduced a representation theoretic lift as a replacement for a microlocal lift.
They sketch, in \cite[Remark 1.7(4) and \S 5.4]{SV07},
a proof that the representation theoretic lift asymptotically
gives the same result as a microlocal lift using \psdiff\ operators.

In this paper, we employ the Riemannian geometric pseudo-differential calculus developed in \cite{Wi,Sh05a,Ha10}.
It has nice equivariance properties.
In particular, a full symbol is invariantly defined, and the symbol and quantization maps are
equivariant under isometries.
It is a useful feature of this quantization, proved in Lemma~\ref{Kap5-lemma-chi-und-symbol}, that the
algebra of invariant fiber-polynomial symbols corresponds to the algebra of invariant differential operators.

We assume the following setting.
Let $X=G/K$ denote a Riemannian symmetric space of noncompact type,
where $G$ is a connected semisimple Lie group with finite center and $K$ a maximal compact subgroup of $G$.
Further, let $\Gamma$ be a co-compact and torsion free discrete subgroup of $G$.
Then we obtain a locally symmetric space $\XG$ as the quotient $\Gamma\backslash X$,
i.e., the double coset space $\Gamma\backslash G/K$.
Let $G=KAN$ be a corresponding Iwasawa decomposition of $G$ and let $M$ denote the centralizer of $A$ in $K$.
The Furstenberg boundary of $X$ can be identified with the flag manifold $B:=K/M$.
Denote by $P=MAN$ the minimal parabolic associated with the Iwasawa decomposition.
Identifying $B$ with $G/P$ we define a $G$-action on $B$.
Under the diagonal action, there is a unique open $G$-orbit $B^{(2)}\cong G/MA$ in $B\times B$.
For rank $1$ spaces $B^{(2)}$ is the set of pairs of distinct boundary points.
In this case each geodesic of $X$ has a unique forward limit point and a unique backward limit point in $B$.
In particular, one can identify $B^{(2)}$ with the space of geodesics.
In higher rank the geometric interpretation is more complicated.
It involves the Weyl chamber flow rather than the geodesic flow.

Joint eigenfunctions come with a spectral parameter $\lambda\in \la^*_\cc$, where $\la$ is the Lie algebra of $A$.
The spectral parameters are unique up to the action of the Weyl group $W$ associated with the Iwasawa decomposition.
The principal part of the spectrum comes from the purely imaginary spectral parameters.
We assume that the spectral parameter of $\phi_h$ is $i\nu_h/h \in i\la^*$, $|\nu_h|=1$.
The Patterson--Sullivan distribution $PS^\Gamma_h\in\Dprime(\Gamma\backslash G/M)$
associated with $\varphi_h$ is constructed as follows.
The Poisson--Helgason transform allows us to write
\begin{equation*}
\varphi_h(x)=\int_B e^{(i\nu_h/h+\rho)A(x,b)} T_h(\intd b), \quad x\in X,
\end{equation*}
where $T_h\in\Dprime(B)$ is the boundary value of $\varphi_h$.
Here, $2\rho$ is the sum of positive restricted roots counted according to multiplicities,
and $A\colon X\times B\to\la$ is the horocycle bracket.
For dealing with non-real $\varphi_h$, it is important that the conjugate of $\varphi_h$
also is a unique transform,
\begin{equation*}
\overline{\varphi_h}(x)=\int_B e^{(-iw_0\cdot\nu_h/h+\rho)A(x,b)} \tilde T_h(\intd b), \quad x\in X,
\end{equation*}
where $\tilde T_h\in\Dprime(B)$.
Here $w_0$ is the longest element of $W$.
The weighted Radon transform $\mathcal R_h:\Ccinfty(G/M)\to \Ccinfty(G/MA)$ is defined by
\[ (\mathcal R_hf)(gMA) =\int_A d_h(gaM,\nu_h)f(gaM)\,\intd a \]
with a weight function related to the the horocycle bracket.
Denote by  $\mathcal R_h':\Dprime(B\times B)\to\Dprime(G/M)$ the dual of $\mathcal R_h$.
The Patterson--Sullivan distribution $PS_h^\Gamma\in \Dprime(\Gamma\backslash G/M)$
is defined as the $\Gamma$-average of ${\mathcal R}_h' (T_h\otimes \tilde T_h)$.

Let $\omega=\lim_h W_h\in\Dprime(T^*\XG)$ be a lifted quantum limit which, after passing to
a subsequence if necessary, has a regular direction $\theta=\lim_h \nu_h$.
In addition, assume
\[ \nu_h=\theta+\bigoh(h)\quad\text{as $h\downarrow 0$.} \]
To link $\omega$ to the sequence $(PS_h^\Gamma)_h$ of Patterson--Sullivan distributions,
we make use of a natural $G$-equivariant map $\Phi\colon G/M\times\la^*\to T^*X$.
For regular $\theta\in \la^*$ this induces a push-forward of distributions,
\[
\Phi(\cdot,\theta)_*:\Dprime(\Gamma\backslash G/M)\to \Dprime(T^* \XG).
\]
Our main result (Theorem~\ref{thm-W-PS-off-diagonal}) can now be stated as follows:
\begin{equation}
\label{W-PS-diagonal}
\omega = \kappa(w_0\cdot\theta) \lim_{h\downarrow 0} (2\pi h)^{\dim N/2} \Phi(\cdot,\theta)_* PS^\Gamma_h
\quad\text{in $\Dprime(T^*\XG)$.}
\end{equation}
Here $\kappa$ is a normalizing function defined in terms of structural data of $X$.
We point out that Theorem~\ref{thm-W-PS-off-diagonal} is more general.
It also describes the situation arising from off-diagonal Wigner distributions
$\big(\Op_{\Gamma,h}(a)\varphi_h\mid \varphi'_h\big)_{L^2(\XG)}$.

If one had a formula interwining Patterson--Sullivan distributions $PS^\Gamma_h$
into lifted Wigner distributions $W_h$, one might be able to deduce \eqref{W-PS-diagonal} as a corollary.
Presumably, an intertwining formula holds only for special \psdiff\ calculi.

The paper is organized as follows.
In Section~\ref{Prelim} we collect various geometric facts needed to construct
the lifted quantum limits and the Patterson-Sullivan distributions.
In particular we discuss the function $\Phi$ and the $G$-orbit $B^{(2)}$.
In Section~\ref{section boundary values} we recall the Helgason-Poisson transform and prove
a regularity theorem of $\Gamma$-invariant boundary values which is instrumental
in proving our main result but also of independent interest (see Theorem~\ref{thm: Polynomial-Abschaetzung}).
In Section~\ref{HHS-Kap4} we give the details of the construction of the Patterson-Sullivan distributions and observe its natural $A$-invariance properties (Remark~\ref{rem:A-equivariance-PS}).
Section~\ref{hhs-oint} provides the technical results on oscillatory integrals which are instrumental in establishing our asymptotic results.
In Section~\ref{hhs-qlim} we describe the lifted quantum limits constructed via
the geometric pseudo-differential calculus and derive their invariance
under the Weyl chamber flow (Theorem~\ref{Kap5-Prop-supp-invar-Wigner}).
In the final Section~\ref{hhs-asymp} we put things together and prove Theorem~\ref{thm-W-PS-off-diagonal}.

%% file: hhs-geom.tex

\section{Geometric Preliminaries}\label{Prelim}

Let $\g$ the Lie algebra of $G$, and $\langle\,,\,\rangle$ the
\emph{Killing form}\index{$\langle\,,\,\rangle$, Killing form} of
$\g$. Let $\theta$\index{$\theta$, Cartan involution} be a
\emph{Cartan involution} of $\g$ such that the form
$(X,Y)\mapsto (X, Y)_\theta:=-\langle X,\theta Y\rangle$ is positive definite on
$\g\times\g$. Let $\g=\lk + \lp$ be the decomposition of $\g$ into
eigenspaces of $\theta$ and $K$ the analytic subgroup of $G$ with
Lie algebra $\lk$. We choose a maximal abelian subspace $\la$ of
$\lp$ and denote by $\la^*$ its dual and $\la^*_{\cc}$ the
complexification of $\la^*$.  Let $A=\exp{\la}$ denote
the corresponding analytic subgroup of $G$ and let $\log$ denote the
inverse of the map $\exp:\la\rightarrow A$.


Given $\lambda\in\la^*$, put $\g_{\lambda} =
\left\{X\in\g \mid (\forall\,H\in\la) [H,X]=\lambda(H)X\right\}$. If
$\lambda\neq0$ and $\g_{\lambda}\neq\left\{0\right\}$, then
$\lambda$ is called a \emph{(restricted) root} and
$m_{\lambda}=\dim(\g_{\lambda})$ is called its \emph{multiplicity}.
Let $\g_{\cc}$ denote the complexification of $\g$ and if
$\mathfrak{s}$ is any subspace of $\g$ let $\mathfrak{s}_{\cc}$
denote the complex subspace of $\g_{\cc}$ spanned by $\mathfrak{s}$.

For $\lambda\in\la^*$ let $H_{\lambda}\in\la$ be determined by
$\lambda(H)=\langle H_{\lambda},H\rangle$ for all $H\in\la$. For
$\lambda,\mu\in\la^*$ we put $\langle\lambda,\mu\rangle:=\langle
H_{\lambda},H_{\mu}\rangle$. Since $\langle\,,\,\rangle$ is positive
definite on $\lp\times\lp$ we set
$|\lambda|:=\langle\lambda,\lambda\rangle^{1/2}$ for
$\lambda\in\la^*$ and $|X|:=\langle X,X\rangle^{1/2}$ for $X\in\lp$.
The $\cc$-bilinear extension of $\langle\,,\,\rangle$ to
$\la_{\cc}^*$ will be denoted by the same symbol.

Let $\la'$ be the open subset of $\la$ where all restricted roots
are $\neq0$. The elements of $\la'$ are called \emph{regular}, and the components of $\la'$ are called \emph{Weyl chambers}. We
fix a Weyl chamber $\la^+$ and call a root $\alpha$ positive ($>0$)
if it is positive on $\la^+$. Let $\la^*_+$ denote the corresponding
Weyl chamber in $\la^*$, that is the preimage of $\la^+$ under the
mapping $\lambda\mapsto H_{\lambda}$. Let $\Sigma$ denote the set of
restricted roots, $\Sigma^+$ the set of positive roots and
$\Sigma^-:= -\Sigma^+$ the set of negative roots.

Let $\Sigma_0=\left\{ \alpha\in\Sigma: \frac{1}{2}\alpha\notin\Sigma
\right\}$ be the set of \emph{indivisible} roots, and put $\Sigma_0^+=\Sigma^+\cap\Sigma_0$,
$\Sigma_0^-=\Sigma^-\cap\Sigma_0$. We set
$\rho:=\frac12\Sigma_{\alpha\in\Sigma^+}m_{\alpha}\alpha$\index{$\rho=2^{-1}\Sigma_{\alpha\in\Sigma^+}m_{\alpha}\alpha$}
and let $N$ denote the analytic subgroup of $G$ with Lie algebra
$\lnn:=\Sigma_{\alpha>0}\g_{\alpha}$. Then
$\overline{\lnn}=\theta(\lnn)=\Sigma_{\alpha<0}\g_{\alpha} $. The
involutive automorphism $\theta$ of $\mathfrak{g}$ extends to an
analytic involutive automorphism of $G$, also denoted by $\theta$,
whose differential at the identity $e\in G$ is the original
$\theta$. It thus makes sense to define $\overline{N}=\theta N$. The
Lie algebra of $\overline{N}$ is $\theta(\lnn)$.

Let $G=KAN$ be the Iwasawa decomposition of $G$ corresponding to the choice of a positive system in $\Sigma$. Writing
\begin{eqnarray}
g=k(g)\exp H(g) n(g),
\end{eqnarray}
where $k(g)\in K$, $H(g)\in\la$, $n(g)\in N$, the functions $k,H,n$ are called the \emph{Iwasawa projections}.
By $M$ we denote the centralizer of $A$ in $K$. Then $P:=MAN$ is a minimal parabolic subgroup of $G$ and $G/P$ is the \emph{Furstenberg boundary} of $X:=G/K$. In view of the Iwasawa decomposition, it can be identified with the flag manifold $B:=K/M$. The group $G$ acts on $G/P$ via $g\cdot xP=gxP$ and $K/M\rightarrow G/P, \, kM\mapsto kP$ is a diffeomorphism (\cite{He01}, p. 407) inverted by $gP\mapsto k(g)M$. Hence this map intertwines the $G$-action on $G/P$ with the action on $K/M$ defined by $g\cdot kM=k(gk)M$. These spaces are thus equivalent for the study of $B=K/M=G/P$.

Let $o:=K\in G/K$ denote the \emph{origin} of the symmetric space $X$ and  $b_+:= M\in K/M$ the canonical base point in $B$.  Then the diagonal action of $G$ on $X\times B = G/K\times G/P  =  G/K\times K/M$ is transitive and the stabilizer of $(o,b_+)$ is $K\cap P=M$, so we can identify $X\times B$ with the space $G/M$ of Weyl chambers as a $G$-space.

Let $M'$ be the normalizer of $A$ in $K$. Then $W:=M'/M$ is the corresponding \emph{Weyl group}. It acts on $\Sigma$ and contains unique element $w_0\in W$ exchanging $\Sigma^+$ and $\Sigma^-$. This element is called the longest element of $W$ and by abuse of notation we will sometimes also denote a representative of $w_0$ in $M'$ by $w_0$. Further, we set $b_-:=w_0\cdot b_+=w_0M\in K/M=B$.

\subsection{The Horocycle Bracket}\label{An equivariance property}

The \emph{horocycle bracket} is defined by
\begin{eqnarray}\label{horocycle bracket}
X\times B \rightarrow \mathfrak{a},\
(gK,kM)\mapsto A(gK,kM):=-H(g^{-1}k).
\end{eqnarray}
Each $(x,b)\in X\times B$ is of the form $(gK,kM)$ and it is easy to see that
\eqref{horocycle bracket} is well-defined. The horocycle bracket is often denoted by
$\left\langle x,b\right\rangle = \left\langle gK,kM\right\rangle = -H(g^{-1}k)$.
In order to avoid confusion with the Killing form we prefer to use the notation $A(x,b)$ over $\langle x,b\rangle$ as in \cite{He94}.
For details on the geometric interpretation of the horocycle bracket we refer to \cite{He94}, Ch. II.

\begin{proposition}
The horocycle bracket  $A\colon X\times B\to\mathfrak{a}$ is invariant under the diagonal action of $K$ on $X\times B$.
\end{proposition}

\begin{lemma}\label{invariance0}
Let $g_1,g_2\in G$, $k\in K$. Then $H(g_1 g_2 k)=H(g_1 k(g_2 k)) + H(g_2 k)$.
\end{lemma}

\begin{proof}
Decompose $g_2k=\tilde{k}\tilde{a}\tilde{n}$ and
$g_1\tilde{k}=k'a'n'$. Then $$H(g_1g_2k) =
H(k'a'n'\tilde{a}\tilde{n}) = H(a'n'\tilde{a}).$$ Since $A$
normalizes $N$ this equals $\log(a') + \log(\tilde{a})$.
\qed
\end{proof}

\begin{lemma}\label{equivariance}
Let $x=hK\in G/K$, $b=kM\in K/M$, $g\in G$. Then
\begin{enumerate}
\item[{\rm(i)}]
$A( g\cdot x,g\cdot b ) = A(x,b) + A(g\cdot o,g\cdot b )$.
\item[{\rm(ii)}]
$A(g^{-1}\cdot o, b)= - A(g\cdot o,g\cdot b)$.
\end{enumerate}
\end{lemma}

\begin{proof}
By definition, $A(g\cdot x,g\cdot b) =
-H(h^{-1}g^{-1}k(gk))$. Then by Lemma \ref{invariance0} applied to
$g_1=h^{-1}g^{-1}$ and $g_2=g$ this equals
\begin{eqnarray*}
-H(h^{-1}g^{-1}gk) + H(gk) = -H(h^{-1}k) + H(gk).
\end{eqnarray*}
For $h=e$ we obtain $A(g\cdot o, g\cdot b) = -H(k) + H(gk) = H(gk)$. Hence
\begin{equation*}
A( g\cdot x,g\cdot b ) - A( g\cdot o,g\cdot b )
=-H(h^{-1}k) = A( hK, kM) = A( x,b) ,
\end{equation*}
which implies (i). For (ii) we use (i) to calculate
$$0=A(o,g\cdot b)=A(g\cdot(g^{-1}\cdot o), g\cdot b)= A(g^{-1}\cdot o, b)+ A(g\cdot o, g\cdot b).$$
\qed
\end{proof}

\begin{lemma}\label{facts}
Let $\gamma, g\in G$. Then
\begin{itemize}
\item[{\rm(i)}] $A( g\cdot o, g\cdot b_+) = H(g)= -A(g\cdot o,b_+)$.
\item[{\rm (ii)}] $A( g\cdot o, g\cdot b_-) = H(gw_0)= -A(g^{-1}\cdot o,b_-)$.
\item[{\rm (iii)}] $H(\gamma g) = H(g) + A(\gamma\cdot o,\gamma g\cdot b_+)$
and
$H(\gamma gw_0) = H(gw_0) + A(\gamma\cdot o,\gamma g\cdot b_-)$.
\end{itemize}
\end{lemma}

\begin{proof}
Parts (i) and (ii) are direct computations. The second part of (iii)
follows from the first part applied to $gw_0$ instead of $g$. For this
assertion, let $z=g\cdot o$. Then by (i)
\begin{eqnarray*}
H(\gamma g) = A( \gamma g \cdot o, \gamma g \cdot b_+)
= A(\gamma\cdot z,\gamma g\cdot b+),
\end{eqnarray*}
which by Lemma~\ref{equivariance} equals
$$A( z,g\cdot b_+) + A( \gamma\cdot o, \gamma g\cdot b_+) = H(g) + A(\gamma\cdot o,\gamma g\cdot b_+).$$
\qed
\end{proof}

\subsection{The Cotangent Bundle and Collective Hamiltonians}\label{subsec:collham}

A detailed study of the cotangent bundle $T^*(X)$ can be found in \cite{Hi05}. We only recall a few facts we will need later on. The $G$-action on $X$ lifts to an action $T(X)$ by taking derivatives and then to an action on $T^*(X)$ by duality.

$T^*(X)$ is a $G$-homogenous vector bundle.
In fact, it  can be written as $G\times_ K\mathfrak{p}^*$, where $K$ acts on $\mathfrak{p}^*$ via the coadjoint representation.
Using the Killing form on $\mathfrak{p}=T_o(X)$, i.e. the invariant Riemannian metric defined by the Killing form, one can identify $T(X)$ and $T^*(X)$. Under this identification adjoint and coadjoint action of $K$ on $\mathfrak{p}$ and $\mathfrak{p}^*$ get identified.

Let $L_g\colon G/K\to G/K$ be the left translation by $g\in G$. Then map
\begin{equation}\label{eq:TX-surjection}
\Phi\colon G/M\times\mathfrak{a}\to T(X)=G\times_K\mathfrak{p}, \  (gM,X)\mapsto dL_g(o)X = [g,X]
\end{equation}
is $G$-equivariant and surjective, but not a covering unless one restricts it to the set $\mathfrak a'$ of regular elements in $\mathfrak{a}$. If one wants to keep $\mathfrak{p}$ and $\mathfrak{p}^*$ apart, the function $\Phi$ is written
\begin{equation}\label{eq:T*X-surjection}
\Phi\colon G/M\times\mathfrak{a}^*\to T^*(X)=G\times_K\mathfrak{p}^*, \  (gM,\theta)\mapsto [g,\theta].
\end{equation}
The fibers of $\Phi$ can be described as follows: $\Phi(gM,\theta)=\Phi(g'M,\theta')$ if and only if there exists a $k\in K$ such that  $g'=gk$ and $k\cdot \theta=\Ad^*(k)\theta= \theta'$. This means
$$\Phi^{-1}([g,\theta])=\{(\tilde gM,\tilde\theta)\in G/M\times \mathfrak{a}^* \mid \exists k\in K : gk=\tilde g, k\cdot \theta=\tilde\theta\}.$$
If $\theta$ is regular, then such a $k$ has to be in $M'$. Therefore, $\tilde gM= gM\cdot w$ and $\tilde \theta =w\cdot \theta$, where $w=kM$ is in the Weyl group $W= M'/M$.

Note that a continuous function $f\colon G/M\times \mathfrak{a}^*\to \cc$ that factors through $\Phi$ will have to satisfy $f(gM\cdot w,w\cdot\theta)=f(gM,\theta)$ for all $w\in W$. But even though the regular elements in $\mathfrak{a}^*$ are dense in $\mathfrak{a}^*$, this condition does not automatically guarantee that $f$ factors through $\Phi$ since the $\Phi$-fibers over the singular points have positive dimension and $W$-invariance cannot guarantee that the function is constant on those fibers as well.

The map $\Phi\colon G/M\times\mathfrak{a}^*\to T^*(X)$ can also be written in terms of the Iwasawa projection (cf. \cite{AS}, \S3.2).

\begin{proposition}\label{prop: Phi-Iwasawa} Consider the function $F\colon X\times B\times \mathfrak a^*\to \rr$ defined by $F(x,b,\theta)=\theta\big(A(x,b)\big)$.
Then $\Phi(x,b,\theta)= dF_x(x,b,\theta)\in T_x^*(X)$.
\end{proposition}

\begin{proof} Identifying $X\times B$ with $G/M$ the map $\Phi$ can be written $\Phi(gM,\theta)=dL_g(o)^{-\top}\theta\in T_{g\cdot o}^*(X)$. Note that  the embedding of $\mathfrak a^*\hookrightarrow \mathfrak p^*$ is given via extension by $0$ on the orthogonal complement of $\mathfrak a$ in $\mathfrak p$. Thus for $v=[x,\xi]=dL_g(o)\xi \in T_x(X)$ we have $\Phi(gM,\theta)(v)=\theta(\xi)$. Therefore it suffices to show that
\begin{equation}\label{AS-equation}
\frac{d}{dt}\Big\vert_{t=0} \theta\Big(A(g\exp t\xi\cdot o, b)\Big)=\theta(\xi)
\end{equation}
for $x=g\cdot o$, $b=g\cdot b_+\in B$,  and $\xi\in\mathfrak p$. To prove this, note first the identity (Lemma~\ref{equivariance})
$$A(g\exp t\xi\cdot o, b)
=A(\exp t\xi\cdot o,b_+)+ A(g\cdot o, b)
= H(\exp t\xi)+ A(g\cdot o, b)
.$$
We claim that
\begin{equation}\label{eq:Heckman}
\lim_{t\to 0} \frac{H(\exp t \xi)}{t} = p_{\mathfrak a}(\xi),
\end{equation}
for all $\xi\in \mathfrak p$, where $p_{\mathfrak a}\colon \mathfrak p \to \mathfrak a$ is the orthogonal projection with respect to the Killing form
(cf. \cite{Heckman},  proof of Theorem 2). Since  $\theta(\xi')=\theta\big(p_{\mathfrak a}(\xi')\big)$ equation \eqref{eq:Heckman} proves \eqref{AS-equation}.

To prove  \eqref{eq:Heckman} it suffices to consider a spanning subset of $\mathfrak p$. If $\xi\in \mathfrak a$, then the claim is clear. If $\xi\in \mathfrak a^\perp $, then we have $\xi=\eta+\theta\eta$ with $\eta\in \mathfrak n$, and one has to show.
$$\lim_{t\to 0} \frac{H(\exp t \xi)}{t} = 0.$$
Writing $\theta\eta+\eta=(\theta\eta-\eta)+2\eta\in \mathfrak k+\mathfrak n$ and using the Campbell--Hausdorff multiplication one calculates
\begin{eqnarray*}
H\big(\exp t(\theta \eta +\eta)\big)
&=&H\big(\exp t\big((\theta\eta-\eta)+2\eta\big)\big)\\
&=&H\big(\exp\big( t(\theta\eta-\eta)*t2\eta+O(t^2)\big)\big)\\
&=&H\big((\exp t(\theta\eta-\eta))(\exp t2\eta) g_t\big)\\
&=& H\big(\exp t(\theta\eta-\eta)\exp t2\eta\big)+ O(t^2)\\
&=&O(t^2),
\end{eqnarray*}
where $g_t$ is a group element differing from the identity by $O(t^2)$.
\qed
\end{proof}

We introduce the involutive algebra of functions on $T^* X$
which are the symbols of invariant differential operators on $X$.

\begin{definition}\label{def:Hi05-material}
Denote by $\mathcal A$ the algebra of $G$-invariant real valued functions in  $C^\infty (T^*X)$ which restrict to polynomials on $\mathfrak p^*=T_o(X)$.
\end{definition}

According to \cite{Hi05}, Theorem 1.1, $\mathcal A$ is finitely generated and its joint level sets are precisely the $G$-orbits in $T^*(X)$. In fact, the proof of that theorem shows that the restriction to $\mathfrak a^*$ induces an isomorphism between $\mathcal A$ and the algebra $I(\mathfrak a^*)$ of Weyl group invariant polynomials on $\mathfrak a^*$ (see also \cite{He00}, Cor.~II.5.12). Note that $\mathcal A$ is also closed under the Poisson bracket $\{f,h\}$.

The \emph{Weyl chamber flow} on $G/M$ is the right $A$-action given by $gM\cdot a:= gaM$. If $X$ is of rank one, i.e. if $\dim_\rr \mathfrak a =1$, it reduces to the geodesic flow on the sphere bundle on $X$.

Given a $G$-invariant function $f\in C^\infty(T^*X)^G$, let $h\in C^\infty(\mathfrak p^*)$ be the restriction to $T_o^*X\cong \mathfrak p^*$. In \cite{Hi05}, \S1, it is shown that the hamiltonian flow $ \rr\times T^*(X)\to T^*(X), (t,\omega)\mapsto \Phi_f^t(\omega)$ associated with $f$ is given by
$$\Phi_f^t([g,\xi]) =\left[ g\exp\big(t\mathrm{grad}\, h(\xi)\big),\xi\right],$$
where the gradient of a function on $\mathfrak p^*$ is taken with respect to the inner product coming from the Killing form. Moreover, considering the restriction of $h$ to $\mathfrak a$ one obtains the following relation between the Weyl chamber flow and the function $\Phi$ from \eqref{eq:T*X-surjection}
\begin{equation}\label{equivA}
\Phi(gM\cdot e^{t\mathrm{grad}\,h(\xi)} ,\xi) =\Phi_f^t\circ\Phi(gM,\xi) \quad \forall \xi\in \mathfrak a^*, gM\in G/M.
\end{equation}
Here it should be noted that $(\mathrm{grad}\,h)|_{\mathrm a^*}=\mathrm{grad}(h|_{\mathrm a^*})$.

 In order to see which Weyl chamber actions $(gM,\xi)\mapsto (gaM,\xi)$ we obtain from \eqref{equivA}, we recall from loc. cit. that
$$\{\mathrm{grad}\,p(\xi)\in \mathfrak a^*\mid p\in I(\mathfrak a^*)\} = \mathfrak a^*$$
if $\xi$ is regular. Note here that the calculations in \cite{Hi05} are done in $T(X)$ rather than $T^*(X)$, but identifying the two bundles via the invariant metric gives the results mentioned above.

We call an element in $(gM,\xi)\in G/M\times \mathfrak a^*$ \emph{regular} if $\xi\in \mathfrak a^*$ is regular. Thus the Hamilton flows associated with functions in $\mathcal A$ preserve the regular elements and produce the entire Weyl chamber flow on the regular elements.

\subsection{Open Cells}

The \emph{Bruhat decomposition} says that $G$ is the disjoint union of the double cosets $PwP$ with $w\in W$, or more precisely, with representatives  of the Weyl group elements in $M'$. Moreover,  $w_0P$ is open in $G$ and this is the only open double coset. In particular $Pw_0P\subseteq G$ is dense. Recall $b_-=w_0M\in K/M$ in $B$ and note that $b_-$ does not depend on the choice of the representative $w_0$ in $M'$.

\begin{proposition}\label{lemma: B2=G/MA} The orbit
$B^{(2)}:=G\cdot (b_+,b_-)$ in $B\times B$ under the diagonal action is open and dense. The stabilizer of $(b_+,b_-)$ is $MA$.
\end{proposition}

\begin{proof}
We claim that
$$G\cdot (b_+,b_-)=\{(h_1P,h_2P)\in B\times B\mid h_2^{-1} h_1\in Pw_0^{-1} P\}.$$
Since $Pw_0P$ is dense and open in $G$ and for $U$ running through a basis of neighborhoods of the identity in $G$, the sets  $h_2^{-1}Uh_1$ form a basis of neighborhoods of $h_2^{-1}h_1$, this set is dense and open in $B\times B$.  Moreover $g\cdot(b_+,b_-)=(b_+,b_-)$ if and only if $g\in P$ and $gw_0\in w_0P$, which is equivalent to $g\in P\cap w_0Pw_0^{-1} =P\cap \theta P= MA$. Thus it only remains to prove the claim. The inclusion ``$\subseteq$'' is clear, so assume that $h_2^{-1}h_1= p_1w_0^{-1}p_2$. Then $h_2=h_1p_2^{-1}w_0p_1^{-1}$ implies $h_2P=h_1p_2^{-1}w_0P$, whence
$$(h_1P,h_2P)=(h_1p_2^{-1}P,h_1p_2^{-1}w_0P)=h_1p_1^{-1}\cdot (b_+,b_-),$$
which proves the claim.
\qed
\end{proof}

\begin{remark}\label{rem:G-Bahnen in B2}
\begin{enumerate}
 \item[(a)]  The Weyl group  $W:=M'/M$ acts from the right on  $G/MA$  via $gMA\cdot wM:= gwMA$ and the induced $W$-action on $B^{(2)}$ is $$(g\cdot b_+,g\cdot b_-)\cdot wM =(gw\cdot b_+, gw\cdot b_-)=(g\cdot (w\cdot b_+), g\cdot (w\cdot b_-)).$$ In particular, we have
     $(b_1,b_2)\cdot w_0M =(b_2,b_1)$
     since $w_0\cdot b_\pm =b_\mp$.

 \item[(b)]  The Weyl group  $W=M'/M$ acts from the right on $G/M$ via $gM\cdot wM:= gwM$ and the induced $W$-action on $X\times B$ is
 $$(g\cdot o,g\cdot b_+)\cdot wM =(gw\cdot o, gw\cdot b_+)=(g\cdot o, k(gw)\cdot b_+).$$

 \item[(c)] $W$ acts also on $G/M$ and $K/M$ from the right such that $K/M\to G/M\to G/MA$ are $W$-equivariant. It is also possible to view the $W$-action on $G/M=X\times B$ as follows: Given $(z,b)\in X\times B\cong G/M$ one finds a corresponding  element $g(z,b)M$ of  $G/M$ and defines $b\cdot_z w= g(z,b)w\cdot b_+$. Then
     $$(z,b)\cdot w =(z, b\cdot_z w),$$
     i.e., the $W$-action on $X\times B$ is a twisted version of the $W$-action on the fibers of  $X\times B\to X$.

\item[(d)]
The argument from the proof of Proposition~\ref{lemma: B2=G/MA} works for any $w\in W$ and proves
$$G\cdot (b_+,w\cdot b_+)=\{(h_1P,h_2P)\in B\times B\mid h_2^{-1} h_1\in Pw^{-1} P\}.$$
Thus the Bruhat decomposition implies that each element $(b,b')\in B^2$ is of the form $g\cdot (b_+,w\cdot b_+)$ for \emph{some} $w\in W$.
     \end{enumerate}
\end{remark}

\begin{remark}\label{rem: section for G/M to G/MA}
It will turn out to be useful to have a smooth section $\sigma\colon G/MA \to G/M$ for the canonical projection $G/M\to G/MA$. To construct $\sigma$ we use the Iwasawa decomposition $G=KNA$ to define a smooth map $\projGM\colon G\to G/M, g=kna\mapsto knM$. Then $kna ma'= km(m^{-1}nm)a'$ for $m\in M$ and $a'\in A $ shows that $\projGM$ factors through the canonical projection $\pi\colon G \to G/MA$. Since $\pi$ is a submersion and $\projGM$ is smooth, the universal property of submersions implies that the resulting map  $\sigma\colon G/MA \to G/M, knaMA\mapsto knM$ is indeed smooth. Using the identifications $G/MA=B^{(2)}$ and  $G/M=X\times B$ from Lemma~\ref{lemma: B2=G/MA} and Remark~\ref{rem:G-Bahnen in B2}, we write
$$\sigma(b,b')=\sigma(g\cdot(b_+,b_-) = kn\cdot (o, b_+) =(kn\cdot o, kn\cdot b_+) =  (z_{b,b'},b), $$
where $(b,b')\mapsto z_{b,b'}$ is defined as the composition of $\sigma$ with the canonical projection $G/M\to G/K$.
\end{remark}

The space $G/M$ can also be interpreted in terms of $B^{(2)}$ as the following proposition shows.

\begin{proposition}\label{rem: Psi}  The map
\begin{eqnarray*}\Psi\colon G/M&\to& B^{(2)}\times A=G/MA\times A \\
kan&\mapsto& (kan\cdot b_+,kan\cdot b_-,a)=(gMA,a)
\end{eqnarray*}
is a diffeomorphism.
\end{proposition}

\begin{proof} Using the properties of the Iwasawa decomposition, it is elementary to check that $\Psi$ is bijective.  Moreover, it is clear that $K\times N\times A\to G/MA\times A,\ (k,n,a)\mapsto (knMA,a)$ is a submersion. So $G\to G/MA\times A,\  g=kan\mapsto (k(ana^{-1})MA, a)$ is a submersion. Thus $\Psi$ is a submersion as well, whence it is a diffeomorphism. \qed
\end{proof}

If we compose $\Psi$ with the canonical embedding $B^{(2)}\times A\hookrightarrow B^2\times A$, we find an embedding $G/M\hookrightarrow B\times B\times A$.

\subsection{Normalization of Measures}\label{Measure theoretic preliminaries}

We briefly recall some normalizations of the measures on the
homogeneous spaces we work with. We follow \cite{He00}. The Killing
form induces Euclidean measures on $A$, $\la$ and $\la^*$. For
$l=\dim(A)$ we multiply these measures by $(2\pi)^{-l/2}$ and obtain
invariant measures $\intd a, \intd H$ and $\intd \lambda$ on $A,\mathfrak{a}$ and
$\mathfrak{a}^*$. This normalization has the advantage that the
Euclidean Fourier transform of $A$ is inverted without a
multiplicative constant. We normalize the Haar measures $dk$ and
$dm$ on the compact groups $K$ and $M$ such that the total measure
is $1$. If $U$ is a Lie group and $L$ a closed subgroup, with left
invariant measures $\intd u$ and $\intd l$, the $U$-invariant measure
$\intd u_L=\intd\,(uL)$ on $U/L$ (if it exists) will be normalized by
\begin{eqnarray}\label{integral quotient}
\int_U f(u)\intd u=\int_{U/L}\left(\int_L f(ul)\intd l\right)\intd u_L.
\end{eqnarray}
This measure exists in particular if $L$ is a compact subgroup of
$U$. In particular, we have a $K$-invariant measure $\intd k_M=\intd\,(kM)$ on
$K/M$ of total measure $1$. We also have a $G$-invariant measure
$\intd x=\intd g_K=d(gK)$ on  $X=G/K$. By uniqueness, $\intd x$ is a constant
multiple of the measure on $X$ induced by the Riemannian structure
on $X$ given by the Killing form. The Haar measures $dn$ and
$\intd\overline{n}$ on the nilpotent groups $N$ and $\overline{N}$ are
normalized such that
\begin{eqnarray}\label{theta(dn)}
\theta(\intd n)=\intd\overline{n}, \hspace{7mm} \int_{\overline{N}}e^{-2\rho(H(\overline{n}))}\intd\overline{n}=1.
\end{eqnarray}
As for $X$ one can also for $N$ consider the Riemannian volume $\intd n_\text{Riem}$ on $N$ given by the left-invariant  Riemannian structure
on $N$ derived from the Killing form. Then $\intd n$ and $\intd n_\text{Riem}$ are proportional and we define the constant $C_N$ via
\begin{equation}\label{eq:CN}
\intd n = C_N\intd n_\text{Riem}.
\end{equation}

\begin{proposition}\label{eta dn =  d oline n}
Set $\eta(n) = w_0 n w_0^{-1}$. Then $\eta(\intd n) = \intd\overline{n}$.
\end{proposition}

\begin{proof}
Since $\eta$ is an automorphism of $G$,  $\eta(\intd n)$ is a Haar measure on $\eta(N)=\theta N= \overline N$. Therefore
$\eta(\intd n) =  c \cdot \intd\overline{n}$ for some constant $c>0$.  We claim that $c=1$. In view of the normalizations \eqref{theta(dn)}
the constant equals
$$
\int_{\eta(N)} e^{-2\rho(H(\eta(n)))} \intd\,(\eta n) = \int_N e^{-2\rho(H(\eta(n)))} \intd n = \int_N e^{-2\rho(H(nw_0^{-1}))} \intd n
$$
and we have
$$ \int_N e^{-2\rho(H(\theta n))} \intd n
=\int_{\theta N} e^{-2\rho(H(\theta n))} \intd\,(\theta n)
=\int_{\overline N} e^{-2\rho(H(\overline n))} \intd\overline{n}
=1.
$$
Let $c_{w_{0}}$ be the conjugation by $w_0$ on $G$. Since $w_0\in K$ and $K$ is the fixed point set of $\theta$, we have
$\theta\circ c_{w_{0}} = c_{w_{0}} \circ \theta$. Thus $\kappa:=\theta\circ c_{w_{0}}$ is an involutive automorphism of $G$, which fixes $N$. This implies $\kappa(dn)= \intd n$, since $\kappa(\intd n) = d \intd n$ with $d>0$ and $d^2=1$. Using
$$\intd n=\kappa(\intd n)= \theta(c_{w_0}(\intd n))$$
we find  $\theta(\intd n) = c_{w_0}(\intd n)$ and calculate
\begin{eqnarray*}
\int_N e^{-2\rho(H(nw_0^{-1}))} \intd n
&=&\int_N e^{-2\rho(H(c_{w_0}n))} \intd n \\
&=&\int_{ c_{w_0}(N)} e^{-2\rho(H(c_{w_0}n))} \intd\,(c_{w_0}n) \\
&=&\int_{\theta N} e^{-2\rho(H(\theta n))} \intd\,(\theta n)\\
&=&1.
\end{eqnarray*}
\qed
\end{proof}

The Haar measure on $G$ can (\cite{He00}, Ch. I, \S5) be normalized such that
\begin{eqnarray}\label{integral formula G}
\int_G f(g)\intd g &=& \int_{KAN}f(kan)e^{2\rho(\log a)}  \intd k \intd a \intd n \\
&=& \int_{NAK}f(nak)e^{-2\rho(\log a)} \intd n \intd a \intd k
\end{eqnarray}
for all $f\in C_c(G)$. Let $f_1\in C_c(AN)$, $f_2\in C_c(G)$, $a\in
A$. Then (\cite{He00}, pp. 182)
\begin{eqnarray} \int_N f_1(na) \intd n = e^{2\rho(\log(a))}\int_N f_1(an) \intd n \end{eqnarray}
and
\begin{eqnarray}\label{follows that} \int_G f_2(g) \intd g = \int_{KNA}f_2(kna) \intd k \intd n \intd a
= \int_{ANK}f_2(ank) \intd a \intd n \intd k.\end{eqnarray}
Let $f_3\in C_c(X)$. It follows from \eqref{follows that} that
\begin{eqnarray}\label{integral formula AN and G/K}
\int_X f_3(x) \intd x = \int_{AN}f_3(an\cdot o) \intd a \intd n.
\end{eqnarray}

%% file: hhs-hebv.tex

\section{Helgason Boundary Values}\label{section boundary values}

\subsection{Eigenfunctions and Poisson Transform}\label{sec:eigenfunctions}

Recall the Harish-Chandra homomorphism $\gamma\colon \mathbb D(X)\to I(\mathfrak a^*)$ which associates a Weyl group invariant polynomial on $\mathfrak a$ with each invariant differential operator on $X= G/K$. The formula
$\chi_\lambda(D)=\gamma(D)(\lambda)$ defines a homomorphism $\chi_\lambda:\mathbb D(X)\to\cc$ for each $\lambda\in \mathfrak a_\cc^*$. In this way one obtains the \emph{joint eigenspace}
\begin{eqnarray*}
\mathcal{E}_{\lambda}(X) = \left\{f\in\mathcal{E}(X)\ \big|\  \big(\forall D\in \mathbb D(X)\big) \
Df=\chi_\lambda(D)f \right\}.
\end{eqnarray*}

Since $\chi_\lambda=\chi_{\lambda'}$ if and only if there exists a $w\in W$ with $\lambda=w\cdot \lambda'$, we see that this is equivalent also to $\mathcal{E}_{\lambda}(X)=\mathcal{E}_{\lambda'}(X)$.

Let $\mathcal{A}(B)$ denote the vector space of analytic functions on $B=K/M$, topologized as in \cite{He94}, \S V.6.1. The \emph{analytic functionals} are (loc. cit.) the functionals in the dual space $\mathcal{A}'(B)$ of $\mathcal{A}(B)$. Fix $\lambda\in\mathfrak{a_\cc}^*$ and recall the set $\Sigma$ of restricted roots. For  $\alpha\in\Sigma$ we write
$\alpha_0:=\alpha/\langle\alpha,\alpha\rangle$. We will need
\emph{Harish-Chandra's} $e$\emph{-functions} (\cite{He94}, p. 163; note that Helgason uses a slightly different notation), defined by
\begin{eqnarray}\label{e function}
e_s^{-1}(\lambda) := \prod_{\alpha\in\Sigma_s^+} \Gamma\left(\frac{m_{\alpha}}{4}+\frac{1}{2}+\frac{\langle \lambda,\alpha_0\rangle}{2}\right) \Gamma\left(\frac{m_{\alpha}}{4}+\frac{m_{2\alpha}}{2}+\frac{\langle \lambda,\alpha_0\rangle}{2}\right),
\end{eqnarray}
where $s\in W$, $\Sigma_s^+:=\Sigma_0^+\cap s^{-1}\cdot\Sigma_0^-$ and
where $\Gamma$ denotes the classical Gamma-function. Note that $\Sigma_{w_0}^+=\Sigma^+_0$ for the longest Weyl group element $w_0$.
Then the fundamental result (\cite{KKMOOT}, see also \cite{Schlichtkrull}, \S 5.4) is:

\begin{theorem}\label{Helgason boundary values}
The Poisson--Helgason transform $P_{\lambda}: \mathcal{A}'(B)\rightarrow\mathcal{E}_{\lambda}(X)$ given by
\begin{eqnarray}\label{Poisson--Helgason transform}
P_{\lambda}(T)(x) := \int_B e^{(\lambda+\rho)A(x,b)} T(\intd b)
\end{eqnarray}
is a bijection if and only if $e_{w_0}(\lambda)\not=0$.
\end{theorem}

Since $\chi_\lambda=\chi_{w\lambda}$ for $w\in W$, one can always assume $\RE \lambda \in \overline{\mathfrak a^*_+}$, so that $e_{w_0}(\lambda)\not=0$. Thus each joint eigenfunction is the Poisson integral of an analytic functional (see \cite{He94}, Theorem V.6.6 and \cite{Schlichtkrull}, Corollary 5.5.4).

One also has a characterization of the class of joint eigenfunctions having distributional boundary values: Let $d_X$ denote the distance function on $X$ and define the space $\mathcal{E}^*(X)$ of smooth functions of \emph{exponential growth} by
\begin{eqnarray}\label{exponential growth}
\mathcal{E}^*(X) := \left\{    f\in\mathcal{E}(X)\mid  (\exists C>0)  \forall x\in X :
|f(x)|\leq C e^{C d_X(o,x)}   \right\}.
\end{eqnarray}
Put $\mathcal{E}_{\lambda}^*(X) := \mathcal{E}^*(X)\cap \mathcal{E}_{\lambda}(X)$. Then one has (cf. \cite{BS87}, Theorem 12.2):

\begin{theorem}\label{distribution}
Suppose that $\lambda\in\mathfrak{a}^*_{\cc}$ is contained in the set
$$\lambdaisom := \left\{\lambda\in\la^*_\cc\ \Big|\  2\frac{\langle \lambda,\alpha\rangle}{\langle \alpha,\alpha\rangle} \not\in -\mathbb N\right\}.$$
Then $P_{\lambda}\colon \mathcal{D}'(B) \to \mathcal{E}_{\lambda}^*(X)$ is a topological isomorphism.
\end{theorem}

For $\lambda\in \lambdaisom$ and $\phi\in \mathcal{E}_{\lambda}^*(X)$ we denote the unique distribution $T\in \mathcal{D}'(B)$ with $P_{\lambda}(T)=\phi$ by $T_{\lambda,\phi}$. We call $T_{\lambda,\phi}$ the \emph{$\lambda$-boundary values}  of $\phi$. Note that $T_{\lambda,\phi}$ actually depends on $\lambda$, since $P_\lambda$ and $P_{\lambda'}$ in general differ even if $\lambda\in  W\cdot \lambda'$.

The space $C^\infty(X)$ has a natural real structure given by the real valued functions. This real structure induces a real structure on the space $\mathbb D(X)$ of invariant differential operators. Here the space $\mathbb D_\rr(X)$ of \emph{real}  invariant differential operators is given as the set of operators in $\mathbb D(X)$ commuting with the complex conjugation on the function spaces. Equivalently, $\mathbb D_\rr(X)$ is the subspace of operators preserving the space of real valued smooth functions.

\begin{proposition}\label{prop:realdiffops}
$\mathbb D(X)$ is spanned by  $\mathbb D_\rr(X)$.
\end{proposition}

\begin{proof}
According to Theorem II.4.9 in \cite{He00} the Harish-Chandra homomorphism maps the algebra  $\mathbb D(X)$  isomorphically onto the algebra  $I(\mathfrak a)$ of $W$-invariant polynomial functions on $\mathfrak a$.  The Harish-Chandra homomorphism is a composition of operations (e.g. taking radial parts) preserving real valued maps (see the arguments leading up to  Theorem II.5.18 in \cite{He00}). Therefore  $\mathbb D_\rr(X)$ gets mapped to the space $I_\rr(\mathfrak a)$ of real valued $W$-invariants. Since $I_\rr(\mathfrak a)$ spans $I(\mathfrak a)$, this implies the claim.
\qed
\end{proof}

Note that the complex conjugation $\overline D$ of $D\in \mathbb D(X)$ is defined by $\overline D(f):= \overline{D(\overline f)}$. Similarly the complex conjugate of a character of $\mathbb D(X)$ is defined by $\overline{\chi}(D):=\overline{\chi(\overline{D})}$.  Therefore,
$D(f)=\chi(D)f$ implies
\begin{equation}\label{eq-Dquer}
D(\overline{f})
=\overline{\overline{D}(f)}
=\overline{\chi(\overline{D})f}
=\overline{\chi(\overline{D})}\,\overline{f}
=\overline{\chi}(D)\overline{f}.
\end{equation}
Since the Harish-Chandra homomorphism commutes with complex conjugation, we have
$$\chi_{\overline\lambda}(D)
=\gamma(D)(\overline\lambda)
=\overline{\overline{\gamma(D)}(\lambda)}
=\overline{\gamma(\overline{D})(\lambda)}
=\overline{\chi_\lambda(\overline{D})}
=\overline{\chi_\lambda}(D).$$
Together, we have proved the first part of the following proposition.

\begin{proposition}\label{rem:realdiffops}
$D(f)=\chi_\lambda(D)f$ implies $D(\overline f)=\chi_{\overline\lambda}(D)\overline f$.
\begin{enumerate}
\item[{\rm(i)}] If $\lambda\in \mathfrak a^*$ is real, then $\mathcal E_\lambda(X)$ is invariant under taking real and imaginary parts. Moreover, $\chi_\lambda(D)$ is real for $D\in \mathbb D_\rr(X)$ and $\mathcal E_\lambda(X)$ is spanned by its real valued elements.
\item[{\rm(ii)}] If $w_0=-\id$, so that  $\gamma(D)(-\lambda) = \gamma(D)(w_0\lambda) = \gamma(D)(\lambda)$, then we have $\overline{\chi_\lambda}=\chi_\lambda$ also for $i\nu=\lambda\in i\mathfrak a^*$. In particular, $\mathcal E_\lambda(X)$ is again invariant under taking real and imaginary parts. Finally, $\chi_\lambda(D)$ is real for $D\in \mathbb D_\rr(X)$ and $\mathcal E_\lambda(X)$ is spanned by its real valued elements.
\item[{\rm(iii)}] Conversely, suppose that there exists a real valued joint eigenvector $\phi\in \mathcal E_\lambda(X)$ with $\lambda\in i\la^*_+$. Then $\lambda$ is contained in the subspace $\ker(w_0+\id)\subseteq \la^*$, which is proper if $w_0\not=-\id$.
\end{enumerate}
\end{proposition}

\begin{proof} To show (ii) we calculate
$$\overline{\chi_\lambda}(D)=\chi_{-\lambda}(D)=\gamma(D)(-\lambda)=\gamma(D)(\lambda)=\chi_\lambda(D).$$
For (iii) we note that $\overline\phi=\phi\in  \mathcal E_\lambda(X)$ implies $\chi_\lambda=\chi_{\overline\lambda}$, whence there exists a $w\in W$ with $-\lambda=\overline\lambda=w\cdot \lambda$.

If $\lambda=i\nu$ is regular, then $\nu$ belongs to an open Weyl chamber in $\la^*$. Since $W$ acts simply transitively on the set of Weyl chambers, we can find a unique $s\in W$ such that  $s\cdot\nu\in \la^*_+$. But then $ sw\cdot \nu = -s\cdot\nu \in -\la^*_+$ so that
$ sws^{-1} (s\cdot \nu)  \in -\la^*_+$. Since  $w_0$ is the unique element in $W$ sending $\la^*_+$ to $-\la^*_+$, this implies
$ sws^{-1} =w_0 $. In particular, if $\nu\in \la^*_+$, i.e. $s=\id$, we find $w=w_0$,  and the claim follows.
\qed
\end{proof}

Recall that complex conjugation on distributions is defined by $\overline{T}(f):=\overline{T(\overline f)}$.

\begin{remark}\label{rem:conj-eigenfunctions}
Let $\lambda\in \lambdaisom$.
Since $\lambdaisom$ is invariant under complex conjugation, also $\overline{\lambda}\in \lambdaisom$.
By Proposition~\ref{rem:realdiffops}, $\phi\in \mathcal{E}_{\lambda}^*(X)$ implies
$\overline\phi\in \mathcal{E}_{\overline\lambda}^*(X)$ and we can write $\phi$ and $\overline\phi$
as Poisson integrals of uniquely determined distributions $T_{\lambda,\phi}$ and $T_{\overline\lambda,\overline\phi}$:
$$
\phi(x)
=P_\lambda(T_{\lambda,\phi})(x)
=\int_B e^{(\lambda+\rho)A(x,b)} T_{\lambda,\phi}(\intd b)
$$
and
\[
\overline \phi(x) = P_{\overline\lambda}(T_{\overline\lambda,\overline\phi})(x)=\int_B e^{(\overline\lambda+\rho)A(x,b)} T_{\overline\lambda,\overline\phi}(\intd b).
\]
On the other hand, taking complex conjugates we find
\begin{eqnarray}\label{eq-Poissonquer}
\overline{\phi}(x)
&=&\overline{T_{\lambda,\phi}\big(e^{(\lambda+\rho)A(x,\cdot)}\big)}
=\overline{T_{\lambda,\phi}}\left(\overline{e^{(\lambda+\rho)A(x,\cdot)}}\right)\\
&=&\int_Be^{( \overline{\lambda}+\rho)A(x,b)} \overline{T_{\lambda,\phi}}(\intd b)
=P_{\overline \lambda}\left(\overline{T_{\lambda,\phi}}\right)(x).\nonumber
\end{eqnarray}
From \eqref{eq-Poissonquer} we deduce $T_{\overline\lambda,\overline\phi}=\overline{T_{\lambda,\phi}}$.
\end{remark}

The following immediate consequence of Remark~\ref{rem:conj-eigenfunctions} will allow us to deal with non-real eigenfunctions (cf. \cite{AZ2}, where a special case is used).

\begin{lemma}\label{lem:Poisson-trafo-quer} Let $\lambda\in\lambdaisom$. If $w\in W$ satisfies $w\cdot\lambda\in\lambdaisom$, then
\[
\overline \phi(x) = P_{w\cdot\overline\lambda}(T_{w\cdot\overline\lambda,\overline\phi})(x)=\int_B e^{(w\cdot \overline\lambda+\rho)A(x,b)} T_{w\cdot\overline\lambda,\overline\phi}(\intd b).
\]
\end{lemma}

\subsection{Spherical Principal Series}

We recall some facts concerning the \emph{principal
series}\index{principal series} representations of $G$. Following
\cite{He94} and \cite{Wil}, let $\nu\in\la^*$ and consider the
representation $\sigma_{\nu}(man)=e^{(i\nu+\rho)\log(a)}$ of
$P=MAN$ on $\cc$. We denote the \emph{induced
representation}\label{induced representation} on $G$ by
$\pi_{\nu}=\Ind_P^G(\sigma_{\nu})$. The \emph{induced
picture}\label{induced picture} of this representation is
constructed as follows: A dense subspace of the representation space
is
\begin{eqnarray*}
H_{\nu}^{\infty} := \left\{ f\in C^{\infty}(G): f(gman)=e^{-(i\nu+\rho)\log(a)}f(g)\right\}
\end{eqnarray*}
with inner product
\begin{eqnarray*}
(f_1\mid f_2) = \int_{K/M}f_1(k)\overline{f_2(k)}\intd k = ({f_1}_{|K}\mid{f_2}_{|K})_{L^2(K/M)}
\end{eqnarray*}
and corresponding norm $\|f\|^2 = \int_{K/M}|f(k)|^2\intd k$. The group
action of $G$ is given by $(\pi_{\nu}(g)f)(x)=f(g^{-1}x)$. The
actual Hilbert space, which we denote by $H_{\nu}$, and the
representation on $H_{\nu}$, which we also denote by
$\pi_{\nu}$, is obtained by completion (cf. \cite{Wil}, Ch. 9).
The representations $\pi_{\nu}$ ($\nu\in\la$) form the
\emph{spherical principal series} of $G$. The representation
$(\pi_{\nu},H_{\nu})$ is a unitary (\cite{He94}, p. 528) and
irreducible (loc. cit. p. 530) Hilbert space representation.

Given $f\in C^{\infty}(K/M)$ we may extend it to a function on $G$
by $\tilde{f}(g)=e^{-(i\nu+\rho)H(g)}f(k(g))$. A direct
computation shows that $\tilde{f}\in H^{\infty}_{\nu}$. On the
other hand, if $f\in H^{\infty}_{\nu}$, then the restriction
$f_{|K}$ of $f$ to $K$ is an element of $C^{\infty}(K/M)$. Moreover,
if $f\in C^{\infty}(K/M)$ and if $\tilde{f}$ is as above, then
$\tilde{f}_{|K}=f$. The mapping $f\mapsto\tilde{f}$ described above
is isometric with respect to the $L^2(K/M)$-norm. We may hence
identify $C^{\infty}(K/M)\cong H_{\nu}^{\infty}$. The advantage
is that the representation space is independent of $\nu$. The
group action on $C^{\infty}(K/M)$ is realized by
\begin{eqnarray}\label{compact model action}
(\pi_{\nu}(g)f)(kM) = f(k(g^{-1}k)M) e^{-(i\nu+\rho)H(g^{-1}k)}.
\end{eqnarray}
This is called the \emph{compact picture} of the (spherical)
principal series. Notice that for $g\in K$ the group action
\eqref{compact model action} simplifies to the left-regular
representation of the compact group $K$ on $K/M$.

Let $\nu\in\la^*$. It follows from
\begin{eqnarray}\label{Poisson follows}
(\pi_{\nu}(g)1)(k)=e^{-(i\nu+\rho)H(g^{-1}k)}=e^{(i\nu+\rho)A(gK,kM)}
\end{eqnarray}
that the Poisson transform $P_{i\nu}(T): G/K\rightarrow \cc$ of
$T\in\mathcal{D}'(B)$ is given by
\begin{eqnarray}\label{Poisson intertwines}
P_{i\nu}(T)(gK) = T(\pi_{\nu}(g)\cdot 1).
\end{eqnarray}

A smooth vector $f\in L^2(K/M)$ is a smooth function on $K/M$.
This follows from the Sobolev lemma, since $f$ and all its derivatives are in $L^2(K/M)$.

\subsection{Regularity of $\Gamma$-invariant Boundary Values}\label{subsec: boundary values}

In this subsection we prove a regularity statement for distribution boundary values of joint eigenfunctions on a compact quotient
$X_{\Gamma}:=\Gamma\backslash X$ of $X$, where $\Gamma$ is a a co-compact, torsion free discrete subgroup  of $G$. Choose a $G$-invariant measure $\nu$ on $\Gamma\backslash G$ such that
\begin{eqnarray}
\int_G f(x) \intd x = \int_{\Gamma\backslash G}\left(\sum_{\gamma}f(\gamma x)\right)\intd \nu(\Gamma x)
\end{eqnarray}
for $f\in C_c(G)$. We will denote the Hilbert space $L^2(\Gamma\backslash G,\nu)$ simply by $L^2(\Gamma\backslash G)$. The $G$-invariance of $\nu$ implies that the equation
\begin{eqnarray*}
(R_{\Gamma}(g)f)(\Gamma x) = f(\Gamma xg)
\end{eqnarray*}
($g,x\in G$, $f\in L^2(\Gamma\backslash G)$) defines a unitary representation $R_{\Gamma}$ of $G$ on $L^2(\Gamma\backslash G)$, which is called the \emph{right-regular representation} of $G$ on $\Gamma\backslash G$.

The action of $G$ on $B$ induces an action on $\mathcal{D}'(B)$ by push-forward: Given $T\in\mathcal{D}'(B)$, a test function $f\in\mathcal{E}(B)$ and $g\in G$, this action is $(gT)(f)=T(f\circ g^{-1})$. When we denote the pairing between distributions and test functions by an integral, we also write $T(\intd\gamma b)$ for $(\gamma T)(\intd b)$.

\begin{remark}\label{rem:Gamma-inv joint eigenfunctions} A joint eigenfunction in $\phi\in L^2(X_\Gamma)$ is automatically smooth, since the Laplace-Beltrami operator is elliptic. Thus we can view it as $\Gamma$-invariant joint eigenfunction $\phi\in \mathcal E_{\lambda}(X)$ which is
automatically contained in $\mathcal E_{\lambda}^*(X)$ since $\Gamma\backslash G$ is compact. According to \cite{He00}, formula (7) in \S IV.5, the eigenvalues of the Laplacian $-\Delta_{X_\Gamma}$ are non-negative and of the form $\langle i\lambda,i\lambda\rangle+|\rho|^2$. Thus, either $\lambda\in i\la^*$ or else $\lambda\in \la^*$ with $|\lambda|\le |\rho|$. In the first case $\lambda$ clearly is contained in $\mathcal A$. In the second case this cannot be guaranteed. The spectral parameters $\lambda$ in $i\la^*$ are called the \emph{principal part} of the spectrum of $L^2(X_\Gamma)$. Thus, for a joint eigenfunction in $\phi\in L^2(X_\Gamma)$ with spectral parameter belonging to the principal part, we have a unique boundary value distribution  $T_{i\nu,\phi}$.
\end{remark}

\begin{proposition}\label{prop:Gamma-inv joint eigenfunctions} Let $\phi\in L^2(X_\Gamma)$ be a joint eigenfunction with spectral parameter $\lambda=i\nu$ belonging to the principal part of the spectrum. Then the boundary value $T_{i\nu,\phi}$ satisfies the invariance condition
\begin{equation}\label{eq:T Gamma Invarianz}
\widetilde{\pi}_{\nu}(\gamma)T_{i\nu,\phi}=T_{i\nu,\phi}\quad \forall \gamma\in\Gamma,
\end{equation}
where $\widetilde{\pi}_{\nu}$ denotes the dual representation on $\mathcal{D}'(B)$ corresponding to the principal series
$\pi_{\nu}$ acting on $H_\lambda^\infty=C^\infty(B)$.

Conversely, if a distribution $T\in \mathcal{D}'(B)$ is invariant under $\widetilde{\pi}_{\nu}(\gamma)$, then
$P_{i\nu}(T)$ is invariant under $\pi_\nu(\gamma)$.
\end{proposition}

\begin{proof}
The equality $\phi(\gamma x)=\phi(x)$ for all $\gamma$ and $x$ implies (recall $A(g\cdot x,g\cdot b ) = A( x,b ) + A( g\cdot o,g\cdot b)$ from Lemma~\ref{equivariance})
\begin{eqnarray*}
\phi(x)
&=& \int_B e^{(i\nu+\rho)A( \gamma\cdot x ,b)} T_{i\nu,\phi}(\intd b)
 = \int_B e^{(i\nu+\rho) A(\gamma\cdot x,\gamma\cdot b)} T_{i\nu,\phi}\big(\intd\,(\gamma\cdot b)\big) \\
&=& \int_B e^{(i\nu+\rho) A( x,b)} e^{(i\nu+\rho)A(\gamma\cdot o,\gamma\cdot b)}
     T_{i\nu,\phi}\big(\intd\, (\gamma\cdot b)\big).
\end{eqnarray*}
By the uniqueness of the boundary value,  we obtain
\begin{eqnarray}\label{boundary values equivariance 1}
T_{i\nu,\phi}\big(\intd\,(\gamma\cdot b)\big) = e^{-(i\nu+\rho)A(\gamma\cdot o,\gamma\cdot b)} T_{i\nu,\phi}(\intd b).
\end{eqnarray}
Now  \eqref{boundary values equivariance 1} and  \eqref{Poisson intertwines} imply the claim.
\qed
\end{proof}

In the situation of Proposition~\ref{rem:Gamma-inv joint eigenfunctions} we consider the mapping
\begin{eqnarray*}
\Phi_{\phi}: H_{\nu}^{\infty} \rightarrow
C^{\infty}(\Gamma\backslash G), \,\,\,\,\,\,\,\, \Phi_{\phi}(f)(\Gamma g)
= T_{i\nu,\phi}(\pi_{\nu}(g)f).
\end{eqnarray*}

\begin{lemma}
$\Phi_{\phi}$ is an isometry w.r.t. the norms of $L^2(K/M)$ and $L^2(\Gamma\backslash G)$.
\end{lemma}

\begin{proof}
The operator $\Phi_{\phi}$ is equivariant with respect to the
actions $\pi_{\nu}$ on $H^{\infty}_{\nu}$ and the right
regular representation of $G$ on $L^2(\Gamma\backslash G)$. We
pull-back the $L^2(\Gamma\backslash G)$ inner product onto the
$(\g,K)$-module $H^{\infty}_{\nu,K}$ of $K$-finite and smooth
vectors (which is dense in $H^{\infty}_{\nu}$, \cite{Wal2}, p.
81):
\begin{eqnarray*}
(f_1\mid f_2)_{2} := ( \Phi_{\phi}(f_1)\mid\Phi_{\phi}(f_2))_{L^2(\Gamma\backslash G)}.
\end{eqnarray*}
Let $f_1\in H^{\infty}_{\nu,K}$. Then $A_{f_1}:
H^{\infty}_{\nu,K} \rightarrow \cc, \, f_2 \mapsto (f_1\mid f_2)_{2}$ is a conjugate-linear, $K$-finite functional on
the $(\g,K)$-module $H^{\infty}_{\nu,K}$. This module is
irreducible and admissible, since $H_{\nu}$ is unitary and
irreducible (\cite{Wal2}, Theorems 3.4.10 and 3.4.11). As $A_{f_1}$
is $K$-finite it is nonzero on at most finitely many $K$-isotypic
components. It follows that there is a linear map
$A:H^{\infty}_{\nu,K}\rightarrow H^{\infty}_{\nu,K}$ such
that for each $f_1\in H^{\infty}_{\nu,K}$ the functional
$A_{f_1}$ equals $f_2\mapsto ( A f_1\mid f_2)_{L^2(K/M)}$.
The equivariance of $\Phi_{\phi}$ and the unitarity of
$\pi_{\nu}$ imply that $A$ is $(\g,K)$-equivariant. Using
Schur's lemma for irreducible $(\g,K)$-modules (\cite{Wal2}, p. 80),
we deduce that $A$ is a constant multiple of the identity and hence
$( \cdot\mid\cdot)_{2}$ is a constant multiple of the
original $L^2(K/M)$-inner product on $H_{\nu,K}^{\infty}$. This
constant is $1$: First, $\Phi_{\phi}(1) = P_{i\nu}(T_{i\nu,\phi}) =
\phi$ is the $K$-invariant lift of $\phi$ to $L^2(\Gamma\backslash
G)$. Then $\|\Phi_{\phi}(1)\|_{L^2(\Gamma\backslash
G)}=1=\|1\|_{L^2(K/M)}$.
\qed
\end{proof}

Let $(y_j)$ and $(x_j)$ be bases for $\lk$ and $\lp$, respectively,
such that $\langle y_j,y_i\rangle = -\delta_{ij}$, $\langle
x_j,x_i\rangle = \delta_{ij}$, where $\langle\,,\,\rangle$ as before denotes
the Killing form. The Casimir operator of $\lk$ is
$\Omega_{\lk}=\sum_i y_i^2$ and the Casimir operator of $\g$ is
\begin{eqnarray*}
\Omega_{\g} = -\sum_{j}x_j^2 + \Omega_{\lk} \in \mathcal{Z}(\g),
\end{eqnarray*}
where $\mathcal{Z}(\g)$ is the center of the universal enveloping algebra $\mathcal{U}(\g)$ of $\g$.

It follows from $T_{i\nu,\phi}(f) = \Phi_{\phi}(f)(\Gamma e)$ that
\begin{eqnarray}\label{first estimate}
|T_{i\nu,\phi}(f)|\leq \|\Phi_{\phi}(f)\|_{\infty}.
\end{eqnarray}
We may now estimate this by a convenient Sobolev norm on
$L^2(\Gamma\backslash G)$. Let $\widetilde{\Delta}$ denote the
Laplace operator of $\Gamma\backslash G$. Then we have
\begin{eqnarray*}
\widetilde{\Delta} = - \Omega_{\mathfrak{g}} + 2\Omega_{\mathfrak{k}}.
\end{eqnarray*}

\begin{definition}
Let $s\in \rr$. The Sobolev space $W^{2,s}(\Gamma\backslash G)$ is
(cf. \cite{Tay81}, p. 22) the space of functions $f$ on
$\Gamma\backslash G$ satisfying $(1+\widetilde{\Delta})^{s/2}(f)\in
L^2(\Gamma\backslash G)$ with norm
\begin{eqnarray*}
\| f \|_{W^{2,s}(\Gamma\backslash G)} = \| (1+\widetilde{\Delta})^{s/2}(f) \|_{L^2(\Gamma\backslash G)}.
\end{eqnarray*}
\end{definition}

Let $m=\dim(\Gamma\backslash G)=\dim(G)$, and let $s>m/2$. The
Sobolev imbedding theorem for the compact space $\Gamma\backslash G$
(\cite{Tay81}, p. 19) states that the identity
$W^{2,s}(\Gamma\backslash G) \hookrightarrow C^0(\Gamma\backslash
G)$ is a continuous inclusion ($C^0(\Gamma\backslash G)$ is equipped
with the usual sup-norm $\|\cdot\|_{\infty}$). It follows that there
exists a $C>0$ such that
\begin{eqnarray}\label{this is 0}
\| \Phi_{\phi}(f)\|_{\infty} \leq C \| \Phi_{\phi}(f)\|_{W^{2,s}(\Gamma\backslash G)} \,\,\,\,\,\,\,\,
\forall f\in C^{\infty}(K/M)  .
\end{eqnarray}

Now we derive the announced regularity estimate for the boundary
values: First, by increasing the Sobolev order, we may assume
$s/2\in\nn$, so
\begin{eqnarray*}
(1+\widetilde{\Delta})^{s/2} =
(1-\Omega_{\mathfrak{g}}+2\Omega_{\mathfrak{k}})^{s/2} \in
\mathcal{U}(\mathfrak{g}).
\end{eqnarray*}
Hence $(1+\widetilde{\Delta})^{s/2}$ commutes with each
$G$-equivariant mapping. Let $f\in H_{\nu}^{\infty}$. Then
\begin{eqnarray}\label{this is 1}
\left\|\Phi_{\phi}(f)\right\|_{W^{2,s}(\Gamma\backslash G)}
&=& \left\|(1+\widetilde{\Delta})^{s/2}\Phi_{\phi}(f)\right\|_{L^2(\Gamma\backslash G)} \nonumber \\
&=& \left\|\Phi_{\phi}((1-\Omega_{\g}+2\Omega_{\lk})^{s/2}(f))\right\|_{L^2(\Gamma\backslash G)} \nonumber \\
&=& \left\| (1-\Omega_{\g}+2\Omega_{\lk})^{s/2}(f) \right\|_{L^2(K/M)}.
\end{eqnarray}
Recall $\pi_{\nu}(\Omega_{\mathfrak{k}})=\Delta_{K/M}$ and
$\Omega_{\g}\in \mathcal{Z}(\g)$. Then \eqref{this is 1} equals
\begin{eqnarray}\label{this is 2}
&&\left\| \sum_{k=0}^{s/2} \binom{s/2}{k} (1+2\Delta_{K/M})^{k}
  (-\Omega_{\g})^{s/2-k}(f)\right\|_{L^2(K/M)} \nonumber \\
&& \hspace{3mm} \leq \, \sum_{k=0}^{s/2} \binom{s/2}{k} \left\|
   (1+2\Delta_{K/M})^{k} (-\Omega_{\g})^{s/2-k}(f)\right\|_{L^2(K/M)}.
\end{eqnarray}
Assume $f\in H_{\nu.K}^{\infty}$ and recall that $\Omega_{\g}$
acts on the irreducible $\mathcal{U}(\g)$-module $H_{\nu,K}^{\infty}$ by
multiplication with the scalar $-(\langle\nu,\nu\rangle +
\langle\rho,\rho\rangle)$ (cf. \cite{Wil}, p. 163), that is
\begin{eqnarray*}
{\Omega_{\g}}_{|H_{\nu,K}^{\infty}} = -\left(\langle\nu,\nu\rangle
+ \langle\rho,\rho\rangle\right) \id_{H_{\nu,K}^{\infty}}.
\end{eqnarray*}
Then \eqref{this is 2} equals
\begin{eqnarray}\label{this is 3}
\sum_{k=0}^{s/2} \binom{s/2}{k} \left\| (1+2\Delta_{K/M})^k
(|\nu|^2+|\rho|^2)^{s/2-k}(f) \right\|_{L^2(K/M)}.
\end{eqnarray}
But $\left(|\nu|^2+|\rho|^2\right)^{-k}\leq 1 + |\rho|^{-s} =:
C'$ ($0\leq k\leq s/2$), so the term in \eqref{this is 3} is bounded
by
\begin{eqnarray}\label{this is 4}
C'\left(|\nu|^2+|\rho|^2\right)^{s/2} \sum_{k=0}^{s/2}
\binom{s/2}{k} \left\| (1+2\Delta_{K/M})^{k}(f) \right\|_{L^2(K/M)}.
\end{eqnarray}
Since $H_{\nu.K}^{\infty}$ is dense in $H_{\nu}^{\infty}$,
this bound holds for all $f\in H_{\nu}^{\infty}$. Using
\eqref{first estimate}-\eqref{this is 4} we get
\begin{eqnarray}\label{estimate}
|T_{i\nu,\phi}(f)| \leq C'\left(|\nu|^2+|\rho|^2\right)^{s/2}
\sum_{k=0}^{s/2} \binom{s/2}{k} \left\| (1+2\Delta_{K/M})^{k}(f)
\right\|_{L^2(K/M)}
\end{eqnarray}
for all $f\in H_{\nu}^{\infty}$ and hence for all $f\in
C^{\infty}(K/M)$.
We set
$$\|f\|_{(s)}:= C'\sum_{k=0}^{s/2} \binom{s/2}{k} \left\| (1+2\Delta_{K/M})^{k}(f)
\right\|_{L^2(K/M)}$$
and note that it is a continuous $C^{\infty}(K/M)$-seminorm independent of $\phi$ and $\nu$.
Since $W$ leaves the norm on $\mathfrak a_\cc^*$ invariant,  \eqref{estimate} yields:

\begin{proposition}\label{Proposition Olbrich}
Let $2s>\dim(G)$ such that $s/2\in\nn$. Then
\begin{eqnarray}
|T_{i\nu,\phi}(f)|\leq (1+|\nu|)^s \|f\|_{(s)} \,\,\,\, \forall \, f\in C^{\infty}(K/M)
\end{eqnarray}
for the distribution boundary values $T_{i\nu,\phi}$ corresponding to a $\Gamma$-invariant joint eigenfunction $\phi\in \mathcal E_{i\nu}(X)$.
\end{proposition}

For $\nu\in\la^*$, let $\mathcal{D}'(B)_{\Gamma}$ denote the space of distributions $T$ on $B$ which satisfy
$\widetilde{\pi}_{\nu}(\gamma)T=T$  for all $\gamma\in\Gamma$. By Proposition~\ref{prop:Gamma-inv joint eigenfunctions},  the Poisson transform $P_{\nu}(T)$ of a distribution $T\in\mathcal{D}'(B)_{\Gamma}$ is a function on the quotient $X_{\Gamma}$. We may hence also define
\begin{eqnarray}\label{my space}
\mathcal{D}'(B)_{\Gamma}^{(1)} := \left\{ T\in\mathcal{D}'(B)_{\Gamma}\ \Big|\  \left\| P_{\nu}(T)\right\|_{L^2(X_{\Gamma})} = 1  \right\}.
\end{eqnarray}
Fix $s$ as in Proposition~\ref{Proposition Olbrich}. Then with
\begin{equation}\label{spectral parameter seminorm}
\mathcal{D}'(B)_{\nu} := \left\{ T\in\mathcal{D}'(B)\ \Big|\  |T(f)|\leq (1+|\nu|)^s \|f\|_{(s)} \,\,\,\, \forall \, f\in C^{\infty}(K/M) \right\}
\end{equation}
the above observations  imply:

\begin{lemma}\label{Lemma Olbrich}
$\mathcal{D}'(B)_{\Gamma}^{(1)}\subseteq\mathcal{D}'(B)_{\nu}$. In other words: There exist $s>0$ and a continuous norm $\|\cdot\|_{(s)}$ on $C^\infty (B\times B)$ such that for any $\Gamma$-invariant joint eigenfunction $\phi\in\mathcal E_{i\nu}(X)$ with spectral parameters $\nu\in\la^*_\cc$ with real part in $\la^*_+$, we have
$$
|T_{i\nu,\phi}(f)| \leq (1+|\nu|)^{s} \|f\|\quad \forall f\in C^{\infty}(B).
$$
The constant $s>0$ and the norm $\|\cdot\|_{(s)}$ are independent of $\phi$ and $\nu$.
\end{lemma}

Each $f\in C^{\infty}(B)\otimes C^{\infty}(B)$ has the form $f =
\sum_{i,j}c_{i,j}f_i\otimes f_j$. We define a cross-norm
$\|\cdot\|$ on $C^{\infty}(B)\otimes C^{\infty}(B)$ by
\begin{eqnarray*}
\|f\| := \inf\Big\{ \sum_{i,j}|c_{i,j}|\|f_i\|_{(s)}\|f_j\|_{(s)} \ \Big|\  f = \sum_{i,j}c_{i,j}f_i\otimes f_j \Big\}.
\end{eqnarray*}
This norm induces a continuous seminorm on the projective tensor
product $C^{\infty}(B)\widehat{\otimes}_{\pi} C^{\infty}(B)$ (cf.
\cite{T}, p. 435). Let $\psi\in \mathcal E_{i\mu}(X)$ denote another $\Gamma$-invariant joint eigenfunction with distribution boundary values
$T_{i\mu,\psi}\in\mathcal{D}'(B)$ and spectral parameter $\mu\in\la^*$.
Given $f = \sum_{i,j}c_{i,j}f_i\otimes f_j \in C^{\infty}(B)\otimes
C^{\infty}(B)$ we obtain
\begin{eqnarray}\label{take inf}
|(T_{i\nu,\phi}\otimes T_{i\mu,\psi})(f)|
&\leq& \sum_{i,j}|c_{i,j}| \cdot |T_{i\nu,\phi}(f_i)| \cdot |T_{i\mu,\psi}(f_j)| \nonumber \\
&\leq& (1+|\nu|)^{s}(1+|\mu|)^{s} \sum_{i,j}|c_{i,j}|\cdot\|f_i\|_{(s)}\cdot\|f_j\|_{(s)},\quad
\end{eqnarray}
which implies (by taking the infimum)
\begin{eqnarray}\label{BtimesB}
|(T_{i\nu,\phi}\otimes T_{i\mu,\psi})(f)| \leq (1+|\nu|)^{s}(1+|\mu)^{s} \|f\|
\end{eqnarray}
for all $f\in C^{\infty}(B)\otimes C^{\infty}(B)$. But
$C^{\infty}(B\times B) \cong C^{\infty}(B) \widehat{\otimes}_{\pi}
C^{\infty}(B)$ (cf. \cite{T}, p. 530) implies that \eqref{BtimesB}
holds for all $f\in C^{\infty}(B\times B)$.

Summarizing we obtain the main result of this section:

\begin{theorem}\label{thm: Polynomial-Abschaetzung} There exist $s>0$ and a continuous norm $\|\cdot\|$ on $C^\infty (B\times B)$ such that for any two $\Gamma$-invariant joint eigenfunctions $\phi\in  \mathcal E_{i\nu}(X)$ and $\psi\in  \mathcal E_{i\mu}(X)$ with spectral parameters $\nu,\mu\in \la^*_\cc$ with real part in $\la^*_+$, we have
$$
|(T_{i\nu,\phi}\otimes T_{i\mu,\psi})(f)| \leq (1+|\nu|)^{s}(1+|\mu)^{s} \|f\|\quad \forall f\in C^{\infty}(B\times B).
$$
The constant $s>0$ and the norm $\|\cdot\|$ are independent of $\phi,\psi,\nu,\mu$.
\end{theorem}

%% file: hhs-psd.tex

\section{Patterson--Sullivan Distributions}\label{HHS-Kap4}

\subsection{Weighted Radon Transforms}

\begin{definition}\label{def:intermediate values}
Given $\nu,\nu'\in\la^*_\cc$, we define $d_{\nu,\nu'}: G/M\rightarrow\cc$ by
\begin{eqnarray}\label{dlambda on G/MA}
d_{\nu,\nu'}(gM) := e^{(i\nu+\rho)H(g)}  e^{(i\nu'+\rho) H(gw_0)}
\end{eqnarray}

\end{definition}

\begin{lemma}\label{equivariance property}
Let $\gamma,g\in G$ and $a\in A$. Then
\begin{enumerate}
\item[{\rm(i)}]
\label{equivariance property2}
$d_{\nu,\nu'}(\gamma gM) = e^{(i\nu+\rho)A(\gamma\cdot o,
\gamma g\cdot b_+)} e^{(i\nu'+\rho) A(\gamma\cdot o, \gamma g\cdot b_-)} d_{\nu,\nu'}(gM).
$\\
\item[{\rm(ii)}]
\label{equivariance property2'}
$d_{\nu,\nu'}(gaM) = e^{i(\nu+w_0\cdot \nu')\log a} d_{\nu,\nu'}(gM).
$
\end{enumerate}
\end{lemma}

\begin{proof}
Part (i) follows from Lemma \ref{facts} and for (ii) we recall that $w_0\cdot\rho =-\rho$ to calculate
\begin{eqnarray*}
d_{\nu,\nu'}(gaM)
&=&e^{(i\nu+\rho)H(ga)}e^{(i\nu'+\rho)H(gaw_0)}\\
&=&e^{(i\nu+\rho)(H(g)+\log a)}e^{(i\nu'+\rho)(H(gw_0)+\log(w_0^{-1}aw_0))}\\
&=&d_{\nu,\nu'}(gM) e^{(i\nu+\rho)\log a}e^{(i\nu'+\rho)\log(w_0^{-1}aw_0)}\\
&=&d_{\nu,\nu'}(gM) e^{i\nu\log a}e^{iw_0\cdot\nu'\log(a)}.
\end{eqnarray*}
\qed
\end{proof}

\begin{definition}
For functions $f$ on $G/M$, the \emph{weighted Radon transform} $\mathcal{R}_{\nu,\nu'}$ on $G/M$ is given by
\begin{eqnarray}\label{repair}
(\mathcal{R}_{\nu,\nu'}f)(g) := \int_A d_{\nu,\nu'}(ga)f(ga)\intd a,
\end{eqnarray}
whenever this integral exists.
\end{definition}

If $\mathcal{R}_{\nu,\nu'}(f)$ exists, then it is a right-$A$-invariant function on $G/M$ and hence a function on $G/MA\cong B^{(2)}$ (cf. Lemma~\ref{lemma: B2=G/MA}).

\begin{lemma}\label{lem:Radon cpt supp}
Let $f\in C_c^{\infty}(G/M)$. Then $\mathcal{R_{\nu,\nu'}}(f)\in C_c^{\infty}(G/MA)= C_c^\infty(B^{(2)})$.
\end{lemma}

\begin{proof}
Projecting the support of $f$ to $G/MA$ we can find a compact subset $C$ of $G/MA$ such that
$$f^a(gM):=f(gaM) = 0$$
for all $a\in A$, whenever $gMA\notin C$. For these $g$ we have $\mathcal{R_{\nu,\nu'}}(f)(g)=0$.
\qed
\end{proof}

\begin{remark}\label{rem:A-equivariance-Radon}\label{A-eigendistributions}
\begin{enumerate}
\item[(i)] \label{rem:Radon extension} Identifying $G/MA$ with $B^{(2)}$ we see that elements of $C^\infty_c(G/MA)$ can be extended by zero to yield elements of $C^\infty_c(B^2)$. In particular we may interpret $\mathcal{R}_{\lambda,\nu'}$ also as an integral transform $C_c^\infty(G/M)\to C^\infty_{c} (B^2)$, where we view $C_c^\infty(B^{(2)})$ as a subset of $C_c^\infty(B^2)$, extending all functions by zero on $B^2\setminus B^{(2)}$.

\item[(ii)]  Lemma~\ref{equivariance property} implies that for $\nu,\nu'\in \mathfrak{a}^*_\cc$ we have
$$\mathcal{R}_{\nu,\nu'}(f^a) = e^{-i(\nu+w_0\cdot\nu')\log a}\mathcal{R}_{\nu,\nu'}(f).$$
In particular, $\mathcal{R}_{\nu,-w_0\cdot \nu}$ is $A$-invariant.

\end{enumerate}
\end{remark}

\begin{proposition}\label{prop: Radon invariance}
Let $\nu,\nu'\in \mathfrak{a}_\cc^*$ and $f\in C_c^\infty(G/M)$. For $\gamma\in G$ set $f_\gamma(gM):=f(\gamma^{-1}gM)$. Then the following equivariance property holds for $(b,b')\in B\times B$.
$$(\mathcal{R}_{\nu,\nu'}f_\gamma)(b,b') =
e^{(i\nu+\rho) A(\gamma\cdot o,b)}e^{(i\nu'+\rho) A(\gamma\cdot o,b')}(\mathcal{R}_{\nu,\nu'}f)(\gamma^{-1}\cdot b,\gamma^{-1}\cdot b').
$$
\end{proposition}

\begin{proof} By Remark~\ref{rem:Radon extension} it suffices to prove the claim for $(b,b') = (g\cdot b_+,g\cdot b_-)$ in $B^{(2)}$, where $gMA$ is determined uniquely by $(b,b')$ (see Proposition~\ref{lemma: B2=G/MA}). Using first Lemma~\ref{equivariance property} and then Lemma~\ref{equivariance} we can calculate
\begin{eqnarray*}
(\mathcal{R}_{\nu,\nu'}f_\gamma)(gMA)
&=& \int_A \, d_{\nu,\nu'}(gaM) \, f(\gamma^{-1}gaM) \, da \\
&=& \int_A d_{\nu,\nu'}(\gamma^{-1}gaM) \, f(\gamma^{-1}gaM) e^{-(i\nu+\rho) A(\gamma^{-1}\cdot o,\gamma^{-1}g\cdot b_+)} \\
&& \hspace{11mm} \times \hspace{2mm} e^{-(i\nu'+\rho) A(\gamma^{-1}\cdot o,\gamma^{-1}g\cdot b_-)} da\\
&=& \int_A d_{\nu,\nu'}(\gamma^{-1}gaM) \, f(\gamma^{-1}gaM) e^{(i\nu+\rho) A(\gamma\cdot o,g\cdot b_+)} \\
&& \hspace{11mm} \times \hspace{2mm} e^{(i\nu'+\rho) A(\gamma\cdot o,g\cdot b_-)} da\\
&=& e^{(i\nu+\rho) A(\gamma\cdot o,g\cdot b_+)}e^{(i\nu'+\rho) A(\gamma\cdot o,g\cdot b_-)}(\mathcal{R}_{\nu,\nu'}f)(\gamma^{-1}gMA).
\end{eqnarray*}
\qed
\end{proof}

If one considers  $\nu,\nu'\in\la^*$, then it is clear from Definition~\ref{def:intermediate values} that $d_{\nu,\nu'}$ as well as its derivatives are of polynomial growth in the spectral parameters. Hence

\begin{proposition}\label{prop: Radon bound}
Let $\chi\in C_c^{\infty}(G/M)$.  For each continuous seminorm $\|\cdot\|_1$ on $C^{\infty}(B^2)$ there is $K>0$ and a continuous seminorm $\|\cdot\|_2$ on $C^{\infty}(G/M)$ such that for all  $f\in C^{\infty}(G/M)$ and all $(\nu,\nu')\in(\la^*)^2$ the estimate
\begin{eqnarray}
\|\mathcal{R}_{\nu,\nu'}(\chi f)\|_1 \leq \left((1+|\nu|)\cdot(1+|\nu'|)\right)^K \|\chi f\|_2
\end{eqnarray}
holds.
\end{proposition}

\subsection{Patterson--Sullivan Distributions on $G/M$ and $\Gamma\backslash G/M$}

\begin{definition}
Fix $\nu,\nu'\in\la_+^*$ and $\phi\in\mathcal{E}^*_{i\nu}(X),\varphi'\in\mathcal{E}^*_{i\nu'}(X)$. Let $T_{i\nu,\phi}$ and $T_{i\nu',\varphi'}$ denote their respective boundary values. The \emph{Patterson-Sullivan distribution} $PS_{\phi,\varphi'}$ on $G/M$ associated to $\phi$ and $\varphi'$ is defined by
\begin{eqnarray}
\langle f,PS_{\phi,\varphi'}\rangle_{G/M} := \int_{B^{(2)}} \, \mathcal{R}_{\nu,-w_0\cdot\nu'}(f)(b,b') \, T_{i\nu,\phi}(\intd b) \,T_{-i w_0\cdot\nu',\overline{\varphi'}}(\intd b'),
\end{eqnarray}
where $f\in C^\infty_c(G/M)$ is a test function. Note that this makes sense since $B^{(2)}$ is open in $B^2$, so the distribution $ T_{i\nu,\phi}(\intd b)\otimes T_{-i w_0\cdot \nu',\overline {\varphi'}}(\intd b')$ on $B^2$ can be restricted to $B^{(2)}$, and $\mathcal{R}_{\nu,-w_0\cdot\nu'}(f)$ is compactly supported in $B^{(2)}$ by Lemma~\ref{lem:Radon cpt supp}.
More precisely, we obtain
\begin{eqnarray}\label{Off-Diagonal Patterson-Sullivan}
\langle f,PS_{\phi,{\varphi'}}\rangle_{G/M} = \int_{B\times B} \,
(\mathcal{R}_{\nu,-w_0\cdot\nu'}f)(b,b') \, T_{i\nu,\phi}(\intd b) \otimes T_{-i w_0\cdot\nu',\overline{\varphi'}}(\intd b').\ \
\end{eqnarray}
\end{definition}

Since boundary values of $\Gamma$-invariant and $L^2(X_{\Gamma})$-normalized eigenfunctions also have polynomial bounds in the eigenvalue parameters, Proposition ~\ref{prop: Radon bound} and Theorem \ref{thm: Polynomial-Abschaetzung} imply the following  estimate:

\begin{proposition}\label{cor:Polynomial Bound} Let $\chi\in C_c^{\infty}(G/M)$. Then there exists $K>0$ and a seminorm $\|\cdot\|$ on $C_c^{\infty}(G/M)$ such that following estimate holds for all $f\in C^{\infty}(G/M)$, all $\nu,\nu'\in\la_+^*$, and all joint eigenfunctions $\phi\in\mathcal{E}^*_{i\nu}(X)$ and ${\varphi'}\in\mathcal{E}^*_{i\nu'}(X)$, which are $\Gamma$-invariant and $L^2(X_{\Gamma})$-normalized:
\begin{equation}\label{define K}
|PS_{\phi,{\varphi'}}(\chi f)| \leq \left((1+|\nu|)\cdot(1+|\nu'|)\right)^K \|\chi f\|.
\end{equation}
\end{proposition}

\begin{proof}
By Theorem \ref{thm: Polynomial-Abschaetzung} and by Proposition
\ref{prop: Radon bound}, we have, for $f\in C^{\infty}(G/M)$,
\begin{eqnarray*}
|PS_{\phi,{\varphi'}}(\chi f)| &=& |(T_{i\nu,\phi}\otimes
T_{- i w_0\cdot\nu',\overline {\varphi'}})(\mathcal{R}_{\nu,-w_0\cdot\nu'}(\chi f))| \\
&\leq& \left((1+|\nu|)\cdot(+|\nu'|)\right)^s
\|\mathcal{R}_{\nu,-w_0\cdot\nu'}(\chi f)\|' \\ &\leq&
\left((1+|\nu|)\cdot(1+|\nu'|)\right)^{s+K} \|\chi f\|, \end{eqnarray*}
where $\|\cdot\|'$ is the fixed seminorm on $C^{\infty}(B\times B)$
constructed in Theorem \ref{thm: Polynomial-Abschaetzung}.  The constants $s$ and $K$ are independent of $f$, since $\|\cdot\|'$ is fixed.
\qed
\end{proof}

The following proposition will allow us to define Patterson--Sullivan distributions also on the quotient $\Gamma\backslash G/M$.

\begin{proposition}\label{off diagonal invariant}
Suppose that $\phi$ and ${\varphi'}$ are $\Gamma$-invariant joint eigenfunctions with spectral parameters $i\nu$ and $i\nu'$ in $i\la^*_+$. Then the distribution $PS_{\phi,{\varphi'}}$ on $G/M$ is $\Gamma$-invariant.
\end{proposition}

\begin{proof}
For $f\in C_c^{\infty}(G/M)$ we calculate, using first  \eqref{boundary values equivariance 1} and then Proposition~\ref{prop: Radon invariance},
\begin{eqnarray*}\label{cancel}
\langle f_{\gamma}, PS_{\phi,{\varphi'}}\rangle_{G/M}
&=& \int_{B\times B} \big(\mathcal{R}_{\nu,-w_0\cdot\nu'}f_\gamma\big)(b,b') \, T_{i\nu,\phi}(\intd b) \otimes T_{-i w_0\cdot\nu',\overline{\varphi'}}(\intd b')\\
&=& \int_{B\times B} \big(\mathcal{R}_{\nu,-w_0\cdot\nu'}f_\gamma\big)(\gamma\cdot(b,b'))
e^{-(i\nu+\rho) A(\gamma\cdot o,\gamma\cdot b)} \nonumber\\
&& \hspace{11mm} \times \hspace{2mm} e^{-(-iw_0\cdot\nu'+\rho) A(\gamma\cdot o,\gamma\cdot b')} \,
T_{i\nu,\phi}(\intd b) \otimes T_{-i w_0\cdot\nu',\overline{\varphi'}}(\intd b')\\
&=& \int_{B\times B} \big(\mathcal{R}_{\nu,-w_0\cdot\nu'}f\big)(b,b')\,T_{i\nu,\phi}(\intd b) \otimes T_{-i w_0\cdot\nu',\overline{\varphi'}}(\intd b')\\
&=&\langle f, PS_{\phi,{\varphi'}}\rangle_{G/M}.
\end{eqnarray*}
\qed
\end{proof}

\begin{remark}\label{rem:A-equivariance-PS} Let $\phi\in\mathcal{E}^*_{i\nu}$ and ${\varphi'}\in\mathcal{E}^*_{i\nu'}$ be $\Gamma$-invariant eigenfunctions. Then by Remark~\ref{rem:A-equivariance-Radon} we see that
\begin{eqnarray}\label{eigendistributions}
\langle f^a, PS_{\phi,{\varphi'}}\rangle_{G/M} = e^{-i(\nu-\nu')\log(a)} \langle f, PS_{\phi,{\varphi'}}\rangle_{G/M}.
\end{eqnarray}
In other words, the $PS_{\phi,{\varphi'}}$ are eigendistributions for the action of $A$ on $G/M$ (given by right-translation). In particular, if $\nu- \nu'=0$, then the associated Patterson--Sullivan distribution is invariant under right-translation by $A$.
\end{remark}

Since $B$ is compact, we can (by using partition of unity) also choose a cutoff $\chi\in C_c^{\infty}(X\times B)$ such that $\sum_{\gamma\in\Gamma}\chi(\gamma\cdot(z,b))=1$. Such a function we call a \emph{smooth fundamental domain cutoff} for $\Gamma$. Let $T\in\mathcal{D}'(X\times B)$ be a $\Gamma$-invariant distribution and $f$ a $\Gamma$-invariant smooth function on $X\times B$. Suppose there is $f_1\in C^\infty_c(X\times B)$ such that $\sum_{\gamma\in\Gamma}f_1(\gamma\cdot(z,b))=f(z,b)$. Then
\begin{eqnarray*}
\langle f_1,T\rangle_{X\times B} &=& \int_{X\times B} \Big\{ \sum_{\gamma\in\Gamma} \chi(\gamma\cdot(z,b)) \Big\} f_1(z,b) \, T(\intd z, \intd b) \\
&=& \int_{X\times B} \sum_{\gamma\in\Gamma} \chi(z,b) \, f_1(\gamma\cdot(z,b)) \, T(\intd z, \intd b).
\end{eqnarray*}
By the invariance of $T$ this equals $\int_{X\times B} \chi(z,b) f(z,b) \, T(\intd z, \intd b)$. We thus have

\begin{proposition}\label{independent}
Let $T\in\mathcal{D}'(G/M)$ be a $\Gamma$-invariant distribution. Let
$f$ be a $\Gamma$-invariant smooth function on $G/M$. Then for any
$f_1,f_2\in C^\infty_c(G/M)$ such that
$\sum_{\gamma\in\Gamma}f_j(\gamma\cdot(z,b))=f(z,b)$ ($j=1,2$) we
have $\langle f_1,T\rangle_{G/M}=\langle f_2,T\rangle_{G/M}$.
\end{proposition}

This proposition implies that the following definition of Patterson--Sullivan distributions on $\Gamma\backslash G/M$ is independent of the choice of a smooth fundamental domain cutoff.

\begin{definition}\label{normalize off-diag} Let $\nu,\nu'\in \la^*_+$.
Suppose that $\phi\in\mathcal{E}^*_{i\nu}(X)$ and ${\varphi'}\in\mathcal{E}^*_{i\nu'}(X)$ are $\Gamma$-invariant joint eigenfunctions. Since $PS_{\phi,{\varphi'}}$ is a $\Gamma$-invariant distribution on $G/M$, the definition descends to the quotient $\Gamma\backslash G/M$ via
\begin{eqnarray}\label{eq:PSGamma}
\langle f,PS^\Gamma_{\phi,{\varphi'}}\rangle_{\Gamma\backslash G/M} := \langle\chi f,PS^\Gamma_{\phi,{\varphi'}}\rangle_{G/M},
\end{eqnarray}
where $\chi$ is a smooth fundamental domain cutoff.



\end{definition}

%% file: hhs-oint.tex

\section{Oscillatory Integrals}\label{hhs-oint}

We deal with the asymptotic behavior of oscillatory integrals
\[\int_X f_h(x,y) e^{i\psi(x,y)/h}\intd x \quad\text{as $h\downarrow 0$.} \]
The parameter $y=(b,b',\nu,\nu')$ ranges in $B^2\times ({\la^*})^2$,
and the phase function arises from non-euclidean plane waves,
\begin{equation}
\label{eq:the-phase-function}
\psi(x,b,b',\nu,\nu') = \nu A(x,b)-(w_0\cdot\nu') A(x,b').
\end{equation}

\subsection{Phase Functions}

We rewrite \eqref{eq:the-phase-function} as follows:
\[ \psi(x,b,b',\nu,\nu') = \nu A(gan\cdot o,g\cdot b_+) - (w_0\cdot\nu') A(gan\cdot o,gw\cdot b_+). \]
Here we used Remark~\ref{rem:G-Bahnen in B2}(d) to write $(b,b')=g\cdot(b_+,w\cdot b_+)$
with $g\in G$ and $w\in W$,
and we defined $a\in A$ and $n\in N$ through $x=gan\cdot o$.
Lemma~\ref{equivariance} and Lemma~\ref{facts} give
\begin{eqnarray*}
A(gan\cdot o,gw\cdot b_+)
&=&A(n\cdot o,w\cdot b_+)+ A(ga\cdot o, gaw\cdot b_+)\\
&=&-H(n^{-1} w)+ H(gaw)\\
&=&- H(n^{-1}w)+ H(gw) + \log(w^{-1}a w)\\
&=&H(gw) + w^{-1}\cdot \log a  -  H(n^{-1}w).
\end{eqnarray*}
In particular, $A(gan\cdot o,g\cdot b_+) = H(ga) = H(g)+ \log a$.
It follows that
\begin{equation}
\label{eq:psi-in-an-coords}
\begin{aligned}
\psi(x,b,b',\nu,\nu') &= \nu H(g)-(w_0\cdot\nu') H(gw) \\
 &\phantom{==} +(\nu-ww_0\cdot\nu')\log a + (w_0\cdot\nu') H(n^{-1}w).
\end{aligned}
\end{equation}
We impose assumptions which will imply that
stationary points of the phase function $\psi$
only arise from the last term. In that context the following set will be important.
\begin{equation}
\label{def-laos}
\laos :=\{(\nu,\nu')\in(\la^*)^2 \mid \forall \1\not=w\in W, \nu\neq w\cdot\nu'\}.
\end{equation}
Notice that $(\nu,\nu)\in\laos$ iff $\nu$ is regular, i.e. $\nu\in\laregs$. Moreover, $(\la_+^*)^2\subseteq \laos$.

We start with a standard observation:

\begin{proposition}  The derivative of the Iwasawa projection $H\colon G\to\la$ is given by
\begin{eqnarray}\label{phase function example}
\intd_{nak}H(nak)(X,Y,Z) &=& \tilde{n}\cdot k^{-1}\cdot a^{-1}\cdot X + \tilde{n}\cdot k^{-1}\cdot Y +  \tilde{n}\cdot Z, \nonumber
\end{eqnarray}
where $nak=\tilde{k}\tilde{a}\tilde{n}\in KAN$.
\end{proposition}

Now we consider the map $\phi^\mu_w$ given by $\phi^\mu_w(n)=\mu\big(H(nw)\big)=\langle H_\mu,H(nw)\rangle$ for $H_\mu\in \mathfrak a$. Then \eqref{phase function example} implies
$$\intd\phi^\mu_w(n)(X)=\langle H_\mu,  \tilde nw^{-1}\cdot X\rangle = \langle w\cdot(\tilde n^{-1}\cdot H_\mu), X\rangle$$
for $X\in \mathfrak n$ and $nw=\tilde k\tilde a\tilde n$.  In order to have a clean description of the critical points of $\phi^\mu_w$ we introduce
$$\Sigma_{w,\pm}:=\{\alpha\in \Sigma^+\mid w\cdot \alpha\in\Sigma^{\pm}\}$$
and set
$N_w:=\exp(\mathfrak n_w)$, where $\mathfrak n_w:=\sum_{\alpha\in \Sigma_{w,+}} \mathfrak g_\alpha$. Note that $N_{w_0}=\{e\}$.

\begin{lemma}\label{prop:N-crit}  For $w\in W$ and $\mu\in \mathfrak a^*_{\mathrm{reg}}$ the set of critical points of the map $\phi^\mu_w\colon N\to\rr$ is $N_w$.
\end{lemma}

\begin{proof}
Writing $\tilde n^{-1}=\exp Y$ we obtain
$$w\cdot(\tilde n^{-1}\cdot H_\mu)=w\cdot H_\mu+ w\cdot (e^{\ad Y}-\id)H_\mu,$$
so that $\intd\phi^\mu_w(n)$ vanishes if and only if the part of $(e^{\ad Y}-\id)H_\mu \in \mathfrak n$ which gets mapped into $\theta\mathfrak n$ by $w$ is zero.

Write $Y=\sum_{\alpha\in \Sigma_{w,+}}Y_\alpha +\sum_{\beta\in \Sigma_{w,-}}Y_\beta$. and let $\beta_0$ be the minimal element $\beta\in\Sigma_{w,-}$ with $Y_{\beta}\not=0$ and note that
$(e^{\ad Y}-\id)H_\mu$ is a finite linear combination of iterated Lie brackets of $Y_\alpha$'s and $Y_\beta$'s. Such an element belong to the root space given by the sum of all the involved $\alpha$'s and $\beta$'s. The minimality condition shows that this root cannot be $\beta_0$. In fact, if it where, no $\beta$'s could occur in the sum, but a sum of roots in $\Sigma_{w,+}$ is again in $\Sigma_{w,+}$. Therefore $(e^{\ad Y}-\id)H_\mu$ contains a summand of the form $-\langle \mu,\beta_0\rangle Y_{\beta_0}$, and if $\langle \mu,\beta_0\rangle\not=0$, then $n$ cannot be a critical point of $\phi_w$. Thus, if $n$ is a critical point, then
$Y\in \mathfrak n_w$ and $\tilde n=\exp(-Y)\in N_w$. This implies $w\tilde n w^{-1}\in N$, and together with $nw=\tilde k\tilde a\tilde n$, also $n =w\tilde n w^{-1}\in N_w\subseteq N\cap wNw^{-1}$, $\tilde a=1$, and $\tilde k= w$.

Conversely, assume that $n\in N_w$. Then $H(nw)=H(w\tilde n)=0$, so that
$$\intd\phi^\mu_w(n)(X)=\langle w\cdot H_\lambda, w\tilde n w^{-1}\cdot X\rangle
=\langle w\cdot H_\mu, n\cdot X\rangle= 0$$
for all $X\in \mathfrak n$, since $n\cdot X\in\mathfrak n$ and $w\cdot H_\mu\in \mathfrak a$.
\qed
\end{proof}

\begin{proposition}[\cite{DKV}]
\label{DKV-phasefcn-on-N}
For $\mu\in\la^*_\mathrm{reg}$ the function
\[
\psi_\mu:N\to\rr, \quad
\psi_\mu(n)= \mu H(n^{-1}w_0),
\]
has $n=e$ as its only critical point.
The Hessian $S(\mu):=\nabla^2\psi_\mu(e)$ is symmetric and non-degenerate.
Its signature and determinant are
\begin{align}
\signature(S(\mu))
   &= \sum_{\alpha\in\Sigma^+} \sign(\langle\mu,\alpha\rangle) \dim(\g_\alpha),\\
\big|\det S(\mu)\big|
   &= \prod_{\alpha\in\Sigma^+} \big|\langle\mu,\alpha\rangle\big|^{\dim(\g_\alpha)}.
\end{align}
\end{proposition}

\begin{proof}
By \cite[Corollary 5.2]{DKV}, the differential of $g\mapsto\mu H(g)$
equals $Y\mapsto \langle Y, n(g)^{-1}\cdot H_\mu\rangle$ at $g\in KAn(g)\subset G$.
A calculation shows that the differential of the embedding $\iota:N\to G$, $n\mapsto n^{-1}w_0$,
is $\intd\iota(n):X\mapsto w_0^{-1} n\cdot (-X)$.
It follows that
\[ \intd \psi_\mu(n): X\mapsto - \langle w_0^{-1} n\cdot X, n(n^{-1} w_0)^{-1}\cdot H_\mu\rangle. \]
In particular, $\intd \psi_\mu(e):X\mapsto - \langle w_0\cdot X, H_\mu\rangle =0$
because $\overline{\lnn}= w_0\cdot\lnn$ is orthogonal to $\la$.
That $e$ is the only critical point of $\psi_\mu$ follows from Lemma~\ref{prop:N-crit}, applied to $w_0\in W$.

By \cite[Lemma 6.1]{DKV}, the Hessian form $\g\times\g\to\rr$ at $g=e$ of $g\mapsto\mu H(g)$ equals
\begin{equation}
\label{H-second-derivative}
(Y,Z)\mapsto \sum_{\alpha\in\Sigma^+} \langle\mu,\alpha\rangle
      \big\langle p_\alpha Y -\theta p_{-\alpha} Y, p_{-\alpha} Z\big\rangle.
\end{equation}
Here $p_\alpha$ is the projection $\g\to \g_\alpha$ corresponding to the direct sum
decomposition $\g=\mathfrak{m}\oplus \la \oplus_{\alpha\in\Sigma} \g_\alpha$.
Composing \eqref{H-second-derivative} with $\intd\iota(e):X\mapsto -w_0\cdot X$,
we deduce
\begin{equation}
\label{psi-mu-Hessian}
\nabla^2\psi_\mu(e)(w_0\cdot X,w_0\cdot Y)
   = \sum_{\alpha\in\Sigma^+} \langle\mu,\alpha\rangle
      \big\langle -\theta p_{-\alpha} X, p_{-\alpha} Y\big\rangle,
\quad X,Y\in \overline{\lnn}.
\end{equation}
By the regularity of $\mu$, we have $\langle\mu,\alpha\rangle\neq 0$.
Since $(X, Y)\mapsto \langle -\theta X, Y\rangle$ is an inner product,
the non-degeneracy of the Hessian and the formulae for the signature
and the determinant are seen after choosing a suitable orthonormal basis
of $\overline{\lnn}=\oplus_{\alpha\in\Sigma^+} \g_{-\alpha}$.
\qed
\end{proof}

\begin{lemma}
\label{geom-no-crit-point}
Assume $(\nu,\nu')\in \laos$.
Then $\intd_x \psi(x,b,b',\nu,\nu') = 0$ iff $\nu' = \nu$,
$(b,b')=g\cdot(b_+,b_-)\in B^{(2)}$, and $x\in gA\cdot o$.
\end{lemma}

\begin{proof}
Suppose $\intd_x \psi(x,b,b',\nu,\nu') = 0$.
Since $\log$ is a diffeomorphism, it follows that $\nu-ww_0\cdot\nu'=0$
in \eqref{eq:psi-in-an-coords}.
Therefore, in view of the assumption, $w=w_0=w_0^{-1}$,
$(b,b')=g\cdot(b_+,w_0\cdot b_+)\in B^{(2)}$, and $\nu=\nu'$.
With these parameters \eqref{eq:psi-in-an-coords} reduces to
\begin{equation}
\label{eq:psi-in-an-coords-w0}
\psi(x,b,b',\nu,\nu') = \nu H(g)-(w_0\cdot\nu') H(gw_0) + (w_0\cdot\nu') H(n^{-1}w_0).
\end{equation}
The remaining assertions follow from Proposition~\ref{DKV-phasefcn-on-N}.
\end{proof}

\subsection{Asymptotics}

It is convenient to have notation for describing asymptotic behavior.
In general, for a given locally convex space $E$, we denote by $h^{-k}E$
the locally convex space of functions $f:I\to E$, $h\mapsto f_h$,
such that $h^k f_h$ is uniformly bounded in $E$.
In particular, $h^0 E$ denotes the space of bounded functions $I\to E$.
Here $I$ is a bounded set of positive reals, having $0$ as a limit point.
The seminorms are $f\mapsto \sup_{h\in I} h^{k}\|f_h\|$,
where $\|\cdot\|$ runs through the seminorms of $E$.
Asymptotic expansions are defined
with respect to the scale $\big(h^{j-k}E\big)_{0\leq j\in \zz}$.
The locally convex space $E=\Ccinfty(\Omega)$ is a regular inductive limit for any second countable smooth manifold $\Omega$.
Therefore, $(f_h)\in h^{-k}\Ccinfty(\Omega)$ iff $(f_h)\in h^{-k}\Ccinfty(K)$
for some compact $K\subset\Omega$.

Lemma~\ref{geom-no-crit-point} and the principle of non-stationary phase
imply the following result.

\begin{lemma}
\label{geom-rapid-decrease}
Let $f_h\in h^0\Ccinfty(X\times B^2\times\laos)$
and compact sets $S\subset X$, $S_B\subseteq B^2$,
such that $S\times S_B$ contains the projections to $X\times B^2$
of the supports of $f_h$.
Assume that $g\cdot(b_+,b_-)\in S_B$ implies $(gA\cdot o)\cap S=\emptyset$.
Then
\begin{equation}
\label{eq:rapid-decrease}
\int_X f_h(x,b,b',\nu,\nu') e^{i\psi(x,b,b',\nu,\nu')/h}\intd x
 \in h^\infty\Ccinfty(B^2\times\laos).
\end{equation}
\end{lemma}

\begin{remark}
\label{no-life-outside-the-diagonal}
Lemma~\ref{geom-no-crit-point} states in particular that the phase function $\psi$
does not have a critical point if $(\nu,\nu')\in\laos$ and $\nu\neq\nu'$.
Therefore, \eqref{eq:rapid-decrease} also holds if the $\laos$-component of
the supports of $f_h$ is contained in a compact subset disjoint to the diagonal.
\end{remark}

We shall be interested in the asymptotic behavior of oscillatory integrals
\begin{equation}
\label{eq:the-oscillatory-integral-off-diagonal}
F_h(b,b',\nu,\nu')=\int_X f_h(x,b,\nu,\nu') e^{\psi(x,b,b',\rho,\rho)}
    e^{i\psi(x,b,b',\nu,\nu')/h}\intd x.
\end{equation}
Lemma~\ref{geom-rapid-decrease} implies that  $F_h(b,b',\nu,\nu')\in h^\infty\Ccinfty(B^2\times\laos)$ if  $f_h\in h^0\Ccinfty(X\times B^2\times\laos)$.

The following construction gives a function
useful for cutting off the integrand near the stationary points.
\begin{lemma}\label{lem:beta}
Let $S\subset X$ compact. There exists $\beta\in \Ccinfty(B^{(2)})\subset\Cinfty(B^2)$
such that $(gA\cdot o)\cap S\neq\emptyset$ implies that $(g\cdot M,g\cdot w_0M)$
is in the interior of the support of $1-\beta$.
Moreover, if we view $\beta\in\Ccinfty(G/MA)$, then the $A$-invariant
lift $\hat\beta\in\Cinfty(G/M)$ of $\beta$ is well-defined.
If $S_A\subset A$ is compact, then the projection of $KS_AN$ to $G/M$
intersects the support of $\hat\beta$ in a compact set.
\end{lemma}

\begin{proof}
In view of the smooth Urysohn lemma, to prove the existence of $\beta$, it suffices to show that
the set of all $gMA\in G/MA\cong B^{(2)}$ for which $gAK/K$ intersects $S$ is compact.
If $S'$ is the preimage of $S$ in $G$ under the canonical projection $G\to G/K$,
then this amounts to the observation that $S'A/MA$ is compact.

An $A$-invariant lift $\hat\beta$ satisfies $\hat\beta(gaM)=\beta(gMA)$.
The existence and uniqueness of $\hat\beta$ is clear.
The support of $\hat\beta$, when viewed as a $MA$-invariant function on $G$,
is contained in $KAS_N$ for some compact $S_N\subset N$.
The assertion about the compactness of the intersection follows.
\qed
\end{proof}

 We introduce a notation for the ordinary \emph{Radon transform}
\begin{equation}
\label{eq:Radon-trafo}
\mathcal R:\Ccinfty(G/M)\to\Ccinfty(G/MA), \quad
\mathcal Rf(gMA)=\int_A f(gaM)\intd a
\end{equation}
and note that, in the situation of Lemma~\ref{lem:beta}, we have $\beta\cdot \mathcal R(f)=\mathcal R(\hat \beta f)$ for all $f\in C_c^\infty(G/M)$.

For $\mu\in\la^*$, we set
\begin{equation}
\label{eq:leading-msp-coeff}
\kappa(\mu)
   = C_N \bigg(\prod_{\alpha\in\Sigma^+} \big|\langle\mu,\alpha\rangle\big|^{\dim(\g_\alpha)}\bigg)^{-1/2}
       e^{\pi i s/4},
\end{equation}
where $C_N$ is defined in \eqref{eq:CN} and the signature
$s = \sum_{\alpha\in\Sigma^+} \sign(\langle\mu,\alpha\rangle) \dim(\g_\alpha)$
is, as a function of $\mu$, constant in each Weyl chamber.

Fix $f_h\in h^0\Ccinfty\big(X\times B\times\laos\big)$ and suppose that $(b,b')=g\cdot(b_+,b_-)=gMA$. Then \eqref{eq:psi-in-an-coords} holds with $w=w_0=w_0^{-1}$, and we have, setting $x=an\cdot o$
\begin{align*}
\psi(g\cdot x,b,b',\nu,\nu')
  & = \nu H(g)-(w_0\cdot\nu')H(gw_0) + (\nu-\nu')\log a+ (w_0\cdot\nu') H(n^{-1}w_0), \\
\psi(g\cdot x,b,b',\rho,\rho) &= \psi(g\cdot x,b,b',\rho) = \rho(H(g)+H(gw_0)) -\rho H(n^{-1}w_0).
\end{align*}
Here we also used $\rho=-w_0\cdot\rho$. Furthermore,
\[
f_h(g\cdot x,b,\nu,\nu') = f_h(gan\cdot o, gan\cdot b_+,\nu,\nu') =f_h(ganM,\nu,\nu').
\]
Using the weight function
\begin{equation}
\label{eq:weight-function-off-diagonal}
d_h(gM,\nu,\nu'):=d_{\nu/h,-w_0\cdot\nu'/h}(gM)= e^{(\frac{i}{h}\nu+\rho)H(g)}  e^{(-\frac{i}{h}w_0\cdot\nu'+\rho) H(gw_0)}
\end{equation}
\eqref{eq:the-oscillatory-integral-off-diagonal}, Lemma~\ref{equivariance property}(ii), and $\intdbar n= e^{-\rho H(n^{-1}w_0)}\intd n$ yield
\begin{eqnarray*}
F_h(b,b',\nu,\nu')
&=&\int_A \int_N f_h(ganM,\nu,\nu')e^{\frac{i}{h}(w_0\cdot\nu')H(n^{-1}w_0) }d_{\nu/h,-w_0\cdot\nu'/h}(gM) \\
&&\qquad\cdot  e^{\frac{i}{h}\big((\nu-\nu')\log a\big)}\intdbar n\intd a\\
&=&\int_A  d_h(gaM,\nu,\nu') \int_N f_h(ganM,\nu,\nu') e^{\frac{i}{h}(w_0\cdot\nu')H(n^{-1}w_0) }\intdbar n \intd a.\\
\end{eqnarray*}

Let $S\subset X$ be a compact set which contains the $X$-projections of the supports of $f_h$. Then consider
\begin{equation}
\label{def-Lh-off-diagonal}
I_h(g,\nu,\nu') := \hat\beta(gM)  \int_N f_h(gnM,\nu,\nu') e^{i(w_0\cdot\nu') H(n^{-1}w_0)/h} \intdbar n,
\end{equation}
where $\beta$ is chosen as in Lemma~\ref{lem:beta}, and  $\hat\beta$ denotes the $A$-invariant lift of $\beta$ to $G/M$.
We have $I_h(g,\nu,\nu')= I_h(gm,\nu,\nu')$ for $m\in M$ since  the weighted measure $\intdbar n$ is  $M$-invariant.
By Lemma~\ref{lem:beta}, Proposition~\ref{DKV-phasefcn-on-N} and the method of stationary phase applied to \eqref{def-Lh-off-diagonal} we get
$I_h\in h^{\dim N/2}\Ccinfty\big(G/M\times\laos\big)$ and an asymptotic expansion
\begin{equation}
\label{Ih-asymp-expansion-off-diagonal}
I_h(gM,\nu,\nu') = \kappa(w_0\cdot\nu')\,(2\pi h)^{\dim N/2} \big(f_h(gM,\nu,\nu') + O(h)\big).
\end{equation}
Here  $\kappa$ is defined by \eqref{eq:leading-msp-coeff}.

The calculation above shows
\begin{eqnarray}\label{eq:betaFh-off-diagonal}
\beta(gMA) F_h(gMA,\nu,\nu')
&=& \int_A d_h(gaM,\nu,\nu') I_h(gaM,\nu,\nu')  \intd a \\
&=&\mathcal R\big(d_hI_h(\cdot,\nu,\nu')\big)(gMA).\nonumber
\end{eqnarray}
On the other hand, Lemma~\ref{geom-rapid-decrease} implies $(1-\beta)F_h\in h^\infty\Cinfty\big(B^2\times\laos\big)$. Together, we obtain
\begin{equation}
\label{Fh-equiv-iintAN-off-diagonal}
F_h - \mathcal R (d_hI_h)\in h^\infty\Ccinfty\big(B^2\times\laos\big).
\end{equation}

We collect these results in the following proposition:

\begin{proposition}
\label{prop:asymptotics-of-Fh-with-cutoff-off-diagonal}
Let $f_h\in h^0\Ccinfty\big(X\times B\times\laos\big)$.
Let $S\subset X$ be a compact set which contains the $X$-projections of the supports of $f_h$.
Choose $\beta$ as in Lemma~\ref{lem:beta}, and denote by $\hat\beta$
the $A$-invariant lift of $\beta$ to $G/M$.
Then $I_h\in h^{\dim N/2}\Ccinfty(G/M\times\laos)$ has the asymptotic expansion
\begin{equation*}
I_h(gM,\nu,\nu') = \kappa(w_0\cdot\nu')\,h^{\dim N/2} \big(f_h(gM,\nu,\nu') + O(h)\big)
\end{equation*}
and the oscillatory integral \eqref{eq:the-oscillatory-integral-off-diagonal} satisfies
\begin{equation*}
F_h - \mathcal R (d_hI_h)\in h^\infty\Ccinfty\big(B^2\times\laos\big).
\end{equation*}
\end{proposition} 

%% file: hhs-qlim.tex

\section{Lifted Quantum Limits}\label{hhs-qlim}

The definition of quantum limits of Wigner measures lifted to the cotangent
bundle, also called semi-classical defect measures,
and the study of their properties is based on semi-classical microlocal analysis.
It is convenient to use a geometric $h$-\psdiff\ calculus.
Refer to \cite{DS}, \cite{EZ03} for $h$-\psdiff\ operators and
to \cite{Sh05a} and \cite[Appendix]{Ha10} for geometric \psdiff\ calculi.
The results in \cite{DS} and \cite{EZ03} are stated for the Weyl quantization.
However, operator classes and principal symbols of operators
do not depend on the chosen quantization.

\subsection{Geometric Pseudo-Differential Calculus}

Let $X$ be a Riemannian manifold.
Denote by $\exp_x:T_x X\to X$ the exponential map of its Levi-Civita connection.
With a symbol $a_h=a(\cdot;h)$ depending on a small parameter $h>0$
we associate a \psdiff\ operator $\Op_h(a_h)$,
\begin{equation}
\label{Op-a-geom}
\Op_h(a_h)u(x)=\int_{T_x^*  X}\int_{T_x X} e^{-i\langle\xi,v\rangle/h} \chi_0(x,v) a_h(x,\xi)
    u(\exp_x v)\intd v\intdbar\xi,
\end{equation}
$x\in X$.
Here $\intdbar\xi =(2\pi h)^{-\dim X}\intd\xi$, and
$\chi_0\in\Cinfty(TX)$ is chosen such that $\chi_0=1$ holds in a neighborhood of the zero section
and that its support is contained in a bounded open neighborhood of the zero section
where the exponential map is injective.
In our applications, the $x$-support of the symbols is compact.

The symbols belong to symbol spaces $S^{m,k}(T^* X)=h^{-k}S^m(T^* X)$.
Often $a_h\in S^{m,k}(T^* X)$ has an asymptotic expansion in powers of $h$,
\[ a_h(x,\xi) \sim \sum_{j\geq 0} h^{-k+j} a_{m-j}(x,\xi), \quad a_\ell\in S^\ell(T^* X). \]
We shall always assume that $a_h$ has a \emph{principal symbol} $h^{-k} a$, i.e.,
$a_h-h^{-k} a\in S^{m,k-1}(T^* X)$ with $a\in S^m(T^* X)$ necessarily uniquely determined.

If $X=\exp_x(B_r)$ is a geodesic ball, then we trivialize the cotangent bundle,
\[
B_r\times T^*_x X\to T^* X, \quad (v,\xi)\mapsto(y, \tau^{T^* X}_{\tofrom{y}{x}}\xi), \;y=\exp_x v.
\]
Here $\tofrom{y}{x}$ denotes the unique geodesic from $x$ to $y$,
and $\tau^{T^* X}_\gamma$ the parallel transport in the cotangent bundle $T^*  X$
along a curve $\gamma$ in $X$.
Using the trivialization, the quantization \eqref{Op-a-geom} is, after a change of variables,
expressed as follows,
\begin{equation}
\label{Op-a-geom-hadamard}
\Op_h(a_h)u(x)=\iint_{T^*  X} e^{-i\langle\xi,\log_x y\rangle/h}
    \psi(x,y) a_h(x,\xi) u(y)\intd y\intdbar\eta.
\end{equation}
Here $\xi=\tau^{T^* X}_{\tofrom{x}{y}}\eta$, $\psi(x,\exp_x v)=\chi_0(x,v)/J(x,v)$,
$\log_x=\exp_x^{-1}$,
and $J(x,v)$ is the determinant of the differential of $\exp_x$ at $v$.
Observe that the phase function in \eqref{Op-a-geom-hadamard} is linear in $\xi$
and generates the conormal bundle of the diagonal in $X\times X$.
Applying \eqref{Op-a-geom-hadamard} with $X$ replaced by convex charts,
one deduces that the definition \eqref{Op-a-geom} leads to known
classes $\Psi^{m,k}(X)$ of $h$-dependent \psdiff\ operators, \cite[Section~8]{EZ03}.
Notice that the cutoff $\chi_0$ in \eqref{Op-a-geom} insures that the operators
$\Op_h(a_h)$ are properly supported.

\begin{remark}\label{rem:diff op}
If the restriction of $a_h$ to each fiber of $T^*X$ is a polynomial,
then $\Op_h(a_h)$ has its Schwartz kernel supported in the diagonal
and thus is a differential operator.
\end{remark}

Modulo residual operators in $\Psi^{-\infty,-\infty}(X)$  the quantization map given
by \eqref{Op-a-geom} is independent of the choice of $\chi_0$.
The symbol isomorphism of the geometric \psdiff\ calculus,
\[
S^{m,k}(T^* X)/S^{-\infty,-\infty}(T^* X) \cong \Psi^{m,k}(X)/\Psi^{-\infty,-\infty}(X),
\]
is given by the quantization map $\Op_h$ and inverted by a symbol homomorphism $\sigma_h$.
On the principal symbol level the rules for compositions and adjoints
agree with those of the Weyl calculus and other quantizations.

The geometric calculus behaves nicely under pullback by isometries.
Let $\varphi:X\to X$ a bijective isometry.
Denote $\varphi^*:\Dprime(X)\to\Dprime(X)$, $u\mapsto u\circ\varphi$,
the pullback operator, and $\varphi^{-*}$ its inverse.
Denote
\[
\intd\varphi^{-\top}:T^* X\to T^* X,
\quad
(x,\intd\varphi(x)^\top\eta)\mapsto(\varphi(x),\eta),
\]
the symplectic map induced by $\varphi$.

\begin{lemma}
\label{lemma-inv-isometries}
For $a_h\in S^{m,k}$,
\begin{equation}
\label{Op-invar}
\varphi^{*}\Op_h(a_h)\varphi^{-*} \equiv \Op_h(a_h\circ \intd\varphi^{-\top})\  \mod \Psi^{-\infty,-\infty}(X).
\end{equation}
Equality holds for differential operators, and if the cutoff in
\eqref{Op-a-geom} satisfies $\chi_0\circ\intd\varphi=\chi_0$.
\end{lemma}
\begin{proof}
Let $a_h\in S^{-\infty,k}$, $u\in\Ccinfty(X)$.
Then \eqref{Op-a-geom} is an absolutely convergent integral.
Set $A=\Op_h(a_h)$ and $B=\Op_h(a_h\circ \intd\varphi^{-\top})$.
Fix $x\in X$, and set $y=\varphi(x)$, $S=\intd\varphi(x)$.
Since $\varphi$ is an isometry, $\varphi(\exp_xv)=\exp_y w$ if $w=S v$.
Using the linear symplectic change of variables $(w,\eta)\mapsto(v,\xi)$,
$w=Sv$ and $\xi=S^\top\eta$, we obtain
\begin{align*}
B\varphi^*  u(x)
  &= \iint_{T_x^*  \times T_x} e^{-i\langle\xi,v\rangle/h} \chi_0(x,v) a_h\big(y,S^{-\top}\xi\big)
          u(\varphi(\exp_x v))\intd v\intdbar\xi  \\
  &= \iint_{T_y^*  \times T_y} e^{-i\langle\eta,w\rangle/h} \chi_1(y,w) a_h(y,\eta)
          u(\exp_y w)\intd w\intdbar\eta.
\end{align*}
Here $\chi_1(y,w)=\chi_0(x,S^{-1}w)$.
Hence $B\varphi^* =\varphi^* A+ \varphi^* R$, where
\[
Ru(y)= \iint_{T_y^*  \times T_y} e^{-i\langle\eta,w\rangle/h} (\chi_1-\chi_0)(y,w) a_h(y,\eta)
          u(\exp_y w)\intd w\intdbar\eta.
\]
Extending by density and continuity to $a\in S^{m,k}$
we obtain $B\varphi^* =\varphi^* A+ \varphi^* R$ with $R\in\Psi^{-\infty,-\infty}(X)$.
Formula \eqref{Op-invar} follows.
Obviously, $R=0$ if $\chi_1=\chi_0$.
To complete the proof we observe that are no non-zero differential
operators in $\Psi^{-\infty,-\infty}(X)$.
\qed
\end{proof}

\subsection{Pseudo-Differential Operators on Locally Symmetric Spaces}

Let $X=G/K$ be a symmetric space of noncompact type as in Section~\ref{Prelim}.
The group $G$ acts on $X$ by left translations which are isometries.
For every $x\in X$, the exponential map $\exp_x:T_xX\to X$ is a diffeomorphism.
Therefore, we define $h$-\psdiff\ operators on $X$ by \eqref{Op-a-geom-hadamard}.

The following lemma relates the geometric \psdiff\ calculus
to Fourier analysis on $X$.
Set $e_{\lambda,b}(x)= e^{(\lambda+\rho)A(x,b)}$
for $x\in X$, $b\in B$, $\lambda\in\lac^*$.
We associate a \emph{non-euclidean symbol} $\tilde{a}_h$ with a symbol $a_h$.
Recall the map $\Phi:(x,b,\theta)\mapsto \intd_x \theta A(x,b)$ from Proposition~\ref{prop: Phi-Iwasawa}.

\begin{lemma}
\label{Kap5-lemma-noneuclid-symb}
Let $a_h\in S^{m,0}(T^*X)$.
Define $\tilde{a}_h$ by
\begin{equation}
\label{noneuclid-symb}
\Op_{h}(a_{h}) e_{i\theta/h,b} = \tilde{a}_h(\cdot,b,\theta) e_{i\theta/h,b}.
\end{equation}
Then $\tilde{a}_h\in h^{0}\Cinfty(X\times B\times \la^*)$.
Moreover, there exists $r_h\in h^{2}\Cinfty(X\times B\times \la^*)$ such that
\begin{equation}
\label{noneuclid-symb-expansion}
\tilde{a}_h(x,b,\theta)= a_h(\xi) + ih (D^{(2)} a_h)(\xi) + r_h(x,b,\theta),
\end{equation}
$\xi=\intd_x \theta A(x,b)\in T_x^* X$.
Here $D^{(2)}$ is a second order differential operator on $T^* X$ with real coefficients.
\end{lemma}
\begin{proof}
Using \eqref{Op-a-geom-hadamard} we write $\tilde{a}_h$
as an oscillatory integral over $(y,\eta)\in T^* X$.
The phase function is
\[
\varphi(x,y,b,\theta,\eta) = -\langle\xi,\log_x y\rangle+\theta(A(y,b)-A(x,b)).
\]
We determine the stationary points of $\varphi$ as a function of $y$ and $\eta$.
First, $\varphi_{\eta}':=\intd_\eta \phi = -\tau^{TX}_{\tofrom{y}{x}}\log_x y$.
It follows that $y=x$ at a stationary point.
Moreover, $\varphi_{\eta\eta}''=0$ and $\varphi_{\eta y}''=-I$ at $y=x$.
Furthermore, $\varphi_{y}'=0$ at $y=x$ implies $\phi_x'(x,b,\theta)=\eta$.
Hence for given $x,b,\theta$ the phase $\varphi$ has the unique stationary
point $(y,\eta)=(x,\phi_x'(x,b,\theta))$ which is non-degenerate.
The signature of the Hessian is zero, and the modulus of its determinant is unity.
Recall the definition of $\intdbar\eta$,
and apply the method of stationary phase.
\qed
\end{proof}

Recall the Definition~\ref{def:Hi05-material} of the algebra $\mathcal A\subset S^\infty(T^* X)$
and the homomorphisms $\chi_\lambda$ from Section~\ref{sec:eigenfunctions}.

\begin{lemma}
\label{Kap5-lemma-chi-und-symbol}
If $p\in \mathcal A$, then $\Op_h(p)\in\DD(X)$.
If $P_h=\Op_h(p_h)\in\DD(X)$, $p_h\in S^{m,0}$,
with principal symbol $p\in \mathcal A$, then
\begin{equation}
\label{Kap5-chi-und-symbol}
\chi_{i\nu/h}(P_h)= p(\nu)+ \bigoh(h) \quad\text{as $h\downarrow 0$,}
\end{equation}
uniformly as $\nu$ stays bounded, $\nu\in\la^*\subset T_o^* X$.
If $P_h^*=P_h$ and $\chi_{i\nu/h}(P_h)$ is real, then \eqref{Kap5-chi-und-symbol} holds with $\bigoh(h)$
replaced by $\bigoh(h^2)$.
\end{lemma}

\begin{proof}
Left translation by an element of $G$ acts as an isometry on $X$.
The first assertion follows from Remark~\ref{rem:diff op}.
Let $P_h=\Op_h(p_h)\in\DD(X)$ with principal symbol $p\in \mathcal A$.
For $\nu\in\la^*$, $h>0$, and $(x,b)\in X\times B$, we have
\[
P_h e_{i\nu/h,b} =\chi_{i\nu/h}(P_h) e_{i\nu/h,b}
 = \tilde{p}_h(x,b,\nu) e_{i\nu/h,b},
\]
where we used \eqref{noneuclid-symb}.
Hence
\[
\chi_{i\nu/h}(P_h) = \tilde{p}_h(x,b,\nu) = p(\nu)+ \bigoh(h),
\]
by \eqref{noneuclid-symb-expansion} and \eqref{AS-equation}.

If $\chi_{i\nu/h}(P_h)$ is real, then $\chi_{i\nu/h}(P_h)= \Re\tilde{p}_h(x,b,\nu) = p(\nu)+ \bigoh(h^2)$. Since $P_h$ is formally self-adjoint, $p$ is real and the subprincipal symbol of $P_h$ is purely imaginary. Therefore, the second term of the stationary phase expansion \eqref{noneuclid-symb-expansion} for $\tilde p_h$ is also purely imaginary. This proves the last assertion.
\qed
\end{proof}

Let $\Gamma$ be a co-compact, torsion-free discrete subgroup of $G$.
The locally symmetric space $\XG=\Gamma\backslash X$ is a Riemannian manifold. 
We denote the quantization map of \eqref{Op-a-geom} by $\Op^{\Gamma}_h$,
if $X$ is replaced by $\XG$.
The notation $\Op_h(a_h)$ continues to denote \psdiff\ operators on $X$.
We identify functions (distributions) on $\XG$ with $\Gamma$-invariant
functions (distributions) on $X$.
Operators on $\Dprime(X)$ which are $\Gamma$-invariant restrict to operators on $\Dprime(\XG)$.
In \eqref{Op-a-geom} the cutoff $\chi_0\in\Cinfty(TX)$ is assumed to
equal unity in a neighbourhood of the zero section.
In addition, we assume that $\chi_0$ is $\Gamma$-invariant, and
is supported in a sufficiently small neighbourhood of the zero section
where the exponential map of $\XG$ is a diffeomorphism. 
By Lemma~\ref{lemma-inv-isometries}, we then have
\begin{equation}
\label{invar-PsDO}
\Op^{\Gamma}_h(a_h)u= \Op_h(a_h)u
\quad\text{for $a_h\in S^{m,k}_\Gamma$, $u\in\Dprime(\XG)$.}
\end{equation}
Here, $S^{m,k}_\Gamma$ denotes the subspace of symbols in
$S^{m,k}\subset\Cinfty(T^* X)$ which are $\Gamma$-invariant.
Denote by $\Psi^{m,k}_\Gamma(X):=\Op^{\Gamma}_h(S^{m,k}_\Gamma)$
the corresponding space of \psdiff\ operators on $\XG$.

Denote by $B(H)$ the algebra, equipped with the operator norm,
of bounded operators on a Hilbert space $H$.
Since $\XG$ is compact, we have $\Psi^{0,0}_\Gamma(X)\subset B(L^2(\XG)$, uniformly bounded in $h$.
This follows from standard $L^2$-continuity properties of \psdiff\ operators.
Moreover, by H\"ormanders's proof \cite[Theorem~18.1.11]{Horm} of $L^2$-continuity
we have, for given $\varepsilon>0$ and uniformly in $h>0$, the estimate
\begin{equation}
\label{Hormander-L2-estimate}
\|\Op^{\Gamma}_h(a_h)\|_{B(L^2(\XG))} \leq (1+\varepsilon)\sup_{T^* X}|a| +\bigoh(\sqrt{h}),
\end{equation}
where $a$ is the principal symbol of $a_h\in S^{0,0}$.
Let $a\in S^0_\Gamma$, $a\geq 0$.
The sharp G{\aa}rding inequality gives that there exists $c>0$ such that
\begin{equation}
\label{sharp-Garding-ineq}
\RE \big(\Op^{\Gamma}_h(a)u\mid u\big)_{L^2(\XG)}\geq -c h\|u\|^2
\end{equation}
for $u\in \Ccinfty(\XG)$.
For a proof see \cite[Theorem~7.12]{DS}, and \cite[Theorem~5.3]{EZ03}.

\subsection{Lifted Quantum Limits}

Every bounded sequence of distributions has a weak*-convergent subsequence.

\begin{lemma}\label{lem: mu probability measure}
Let $(\varphi_j)_j, (\phi'_j)_j$ be bounded sequences in $L^2(\XG)$, $0<h_j\to 0$.
Set
\[
W_j(a)=\big(\Op^{\Gamma}_{h_j}(a)\varphi_j\mid \phi'_j\big)_{L^2(\XG)},
\quad a\in\Ccinfty(T^* \XG).
\]
Then $(W_j)_j$ is a bounded sequence in $\Dprime(T^* \XG)$.
Assume that $\omega=\lim_j W_j$ in $\Dprime(T^* \XG)$ as $j\to\infty$.
Then $\omega$ is a Radon measure on $T^* \XG$ of finite total variation, and
\begin{equation}
\label{mu-Wigner}
\int_{T^* \XG} a\intd \omega
   = \lim_{j\to\infty} \big(\Op^{\Gamma}_{h_j}(a_{h_j})\varphi_j\mid \phi'_j\big)_{L^2(\XG)}
\end{equation}
if $a_{h_j}\in S^{0,0}_\Gamma$ has principal symbol $a\in S^0_\Gamma$.
If $\|\varphi_j\|_{L^2(\XG)}=1$ and $\phi'_j=\varphi_j$, then $\omega$ is a probability measure.
\end{lemma}

\begin{proof}
Since $\Op^{\Gamma}_h$ maps $S^{0,0}_\Gamma$ continuously into $B(L^2(\XG))$,
the boundedness of $(W_j)_j$ follows.
Now assume $\lim_{j\to\infty} W_j=\omega$.
Let $M\geq\sup_j(\|\varphi_j\|, \|\phi'_j\|)$.
It follows from \eqref{Hormander-L2-estimate} that
$\limsup_j|W_j(a)|\leq M^2 \sup_{T^*  X} |a|$,
implying that $\omega$ is a Radon measure of total variation $\leq M^2$.
Now assume $\|\varphi_j\|_{L^2(\XG)}=1$ and $\phi'_j=\varphi_j$.
Thus, we can choose $M=1$.
If $0\leq a\in\Ccinfty(T^* \XG)$, then it
follows from the sharp G{\aa}rding inequality \eqref{sharp-Garding-ineq}
and $\IM \big(\Op_h(a)u\mid u\big)=\bigoh(h)$ that $\omega(a)\geq 0$.
Notice $\omega(1)=1$.
Thus $\omega$ is a probability measure.
\qed
\end{proof}

Let $(\varphi_j)_j\subset L^2(\XG)$ be a sequence of normalized
joint eigenfunctions of the algebra $\DD(X)$ of invariant differential operators on $X$
with associated spectral parameters $\lambda_j\in\lac^*$,
$D\varphi_j=\chi_{\lambda_j}(D)\varphi_j$ if $D\in\DD(X)$.
Let $\Delta=\Delta_{\XG}\in\DD(X)$ denote the Laplacian on $\XG$.
The eigenvalues
$\chi_{\lambda_j}(-\Delta)=-\langle\lambda_j,\lambda_j\rangle+|\rho|^2$
are non-negative.
We restrict attention to the principal spectrum, \cite{DKV1}.
Therefore, assume that $\lambda_j\in i\la^*$.
Set $\lambda_j=i\nu_j/h_j$ with unit vectors $\nu_j\in\la^*$, $h_j=|\lambda_j|^{-1}$.
We say that $(\varphi_j)_j$ has \emph{lifted quantum limit}
$\omega$ if the sequence of distributions $W_j\in\Dprime(T^* \XG)$,
\[
W_j(a)=\big(\Op^{\Gamma}_{h_j}(a)\varphi_j\mid \varphi_j\big)_{L^2(\XG)},
\quad a\in\Ccinfty(T^* \XG),
\]
converges, $\omega=\lim_{j\to\infty} W_j$.
Passing to a subsequence, we can assume that $\theta=\lim_{j\to\infty} \nu_j\in\la^*_\cc$ exists.
Following \cite{AS} we then say that $\omega$ is the lifted quantum limit in the direction $\theta$.
The distributions $W_j$ are lifts of the Wigner measures $w_j=|\varphi_j|^2\intd x$
under the canonical projection $\pi:T^* \XG\to\XG$, $\pi_* W_j=w_j$.

In addition, we assume that the sequence $(h_j)_j$ is strictly decreasing.
We can then use $h$ as a subscript, removing $j$ from the notation.
In particular, we denote the spectral parameters $i\nu_h/h$,
and we write
\begin{equation}
\label{def-wigner-measure}
\int_{T^*\XG} a\intd\omega =\lim_{h\downarrow 0}\big(\Op^{\Gamma}_h(a)\varphi_h\mid \varphi_h\big)_{L^2(\XG)}.
\end{equation}

Using the metric tensor we regard the unit sphere bundle $S^*\XG$
as a subset of the cotangent bundle $T^*\XG$.
Then, in view of the results recalled in Subsection~\ref{subsec:collham},
propagation of singularities and Lemma~\ref{Kap5-lemma-chi-und-symbol}
allow us to prove the following invariance properties of lifted quantum limits.

\begin{theorem}[{\cite[Theorem 1.6(3)]{SV07}, \cite[Theorem 1.3]{AS}}]
\label{Kap5-Prop-supp-invar-Wigner}
Assume that $(\varphi_h)_h$ has the lifted quantum limit $\omega$. Then $\supp(\omega)\subset S^*\XG$, and $\omega$ is invariant under the geodesic flow.
Moreover, $\supp \omega$ is contained in a joint level set of $\mathcal A$, i.e. in a $G$-orbit in $S^*\XG$. Moreover, for every $p\in \mathcal A$, $\omega$ is invariant under the Hamilton flow generated by $p$. If the direction $\theta\in\la^*$ of $\omega$ is regular, then $\omega$ is $A$-invariant.
\end{theorem}

\begin{proof}
We can assume that $\omega$ is a lifted quantum limit in the direction
$\theta=\lim_{h\downarrow 0} \nu_h\in\la^*\subset T_o^* X$.

Note that $-h^2\Delta_X=\Op_h(g)$, where $g\in\mathcal A$ is the metric form,
$g(\xi)=|\xi|^2$, $\xi\in T^* X$.
Since $|\nu_h|=1$ and $\nu_h$ is real,
$\chi_{i\nu_h/h}(-h^2\Delta) = 1 +h^2|\rho|^2$.
Hence
\[
\|h^2\Delta_{\XG}\varphi_h +\varphi_h\|_{L^2} = \bigoh(h^2).
\]
It follows from \cite[Theorem~5.4]{EZ03} that the support of $\omega$ is contained
in $S^* \XG=g^{-1}(1)$ because this is the characteristic variety
of the $h$-differential operator $h^2\Delta_{\XG} +1$.
The invariance under the geodesic flow follows from \cite[Theorem~5.5]{EZ03}.%

Let $p\in\mathcal A$.
Set $P_h=\Op_h(p)\in\DD(X)$.
Choose an integer $m$ such that the order of $P_h$ is $<2m$.
Define the $h$-differential operator
$Q_h=\Op_h(p-p(\xi)+g^m-1)\in\DD(X)$, $0<h<1$.
By Lemma~\ref{Kap5-lemma-chi-und-symbol} we have
\[
\|Q_{h} \varphi_h\|_{L^2}= \bigoh(h)\quad\text{as $h\to 0$.}
\]
It follows from \cite[Theorem~5.4]{EZ03} that $\supp(\omega)$ is contained
in the characteristic variety of $Q_h$. The latter intersected with
the unit sphere bundle is contained in the level set $p^{-1}(p(\theta))$.
This proves that $\supp(\omega)$ is contained in a joint level set.

We prove that $\omega$ is invariant under the Hamilton flow generated by $p$.
Adding a constant to $p$ if necessary, we may assume that $p=1$ on $\supp(\omega)$.
By selfadjointness, the eigenvalues of $P_{h}^* P_{h} +(-h^2\Delta)^{2m}$ are real so by the last assertion of  Lemma~\ref{Kap5-lemma-chi-und-symbol} we have
\[
\|\big(P_{h}^* P_{h} +(-h^2\Delta)^{2m} - 2\big)\varphi_h\|_{L^2}=\bigoh(h^2).
\]
By \cite[Theorem~5.5]{EZ03}
we have, for every $a\in\Ccinfty(T^* X)$,
\[
0=\int \{p^2 +g^{2m},a\} \intd\omega = 2\int \{p,a\} \intd\omega.
\]
Here we used the invariance of $\omega$ under the geodesic flow and $p=1$ on $\supp\omega$.
This proves the invariance of $\omega$ under the Hamilton flow generated by $p$.%

Recall from Subsection~\ref{subsec:collham}
that \eqref{equivA} intertwines the Weyl chamber flow with certain Hamilton flows.
Indeed the last statement of that subsection says that because $\theta$ is regular,
each one-parameter subgroup of the Weyl chamber flow can be realized as a Hamilton flow associated with a function in $\mathcal A$. Thus, the $A$-invariance follows.
\qed
\end{proof}

%% file: hhs-asymp.tex

\section{Spectral Directions and Asymptotics}\label{hhs-asymp}

We study the asymptotic behavior of the principal spectrum  of $\XG$ which correspnds to spectral parameters $\lambda\in i\la^*$; see \cite{DKV1}.
Let $(\varphi_h)_h, (\varphi'_h)_h\subset L^2(\XG)$ be sequences of normalized joint eigenfunctions,
with purely imaginary spectral parameters $i\nu_h/h\in i\la^*$, $i\nu'_h/h\in i\la^*$.
The Poisson--Helgason transform~\eqref{Poisson--Helgason transform} gives unique representations,
\begin{equation}\label{eq:hi_h-off-diagonal}
\varphi_h(x)=\int_B e^{(i\nu_h/h+\rho)A(x,b)} T_{i\nu_h/h,\phi_h}(\intd b), \quad x\in X.
\end{equation}
We use Lemma~\ref{lem:Poisson-trafo-quer}
to pick a suitable representation of $\overline{\varphi'_h}$ as a Poisson integral:
\begin{equation}\label{eq:hi_h-quer}
\overline{\varphi'_h}(x)=\int_B e^{(-iw_0\cdot\nu'_h/h+\rho)A(x,b')} T_{-iw_0\cdot \nu'_h/h,\overline{\varphi'_h}}(\intd b'), \quad x\in X.
\end{equation}
This reduces to \eqref{eq:hi_h-off-diagonal} if $w_0=-\id$ and if $\phi_h=\varphi'_h$ is real valued.
To simplify our notation we write $T_h$ and $\tilde T_h$ for  $T_{i\nu_h/h,\phi_h}$
and $T_{-iw_0\cdot \nu'_h/h,\overline{\varphi'_h}}$, respectively.

\begin{lemma}\label{lem:direct-calculation-off-diagonal}
Let $\chi\in\Ccinfty(X)$ real-valued, and $a_h\in S^{0,0}(T^*X)$.
Then
\begin{equation}
\label{direct-calculation-off-diagonal}
\big(\Op_{h}(a_h) \varphi_h\mid \chi\varphi'_h\big)_{L^2(X)}
  = \int_{B^2} F_h(b,b',\nu_h,\nu'_h) \,  T_h(\intd b)\otimes \tilde T_h(\intd b'),
\end{equation}
where
\begin{equation}
\label{direct-calculation-Fh-off-diagonal}
F_h(b,b',\nu,\nu') = \int_X \chi(x) \tilde{a}_h(x,b,\nu) e^{\psi(x,b,b',\rho,\rho)}
  e^{i\psi(x,b,b',\nu,\nu')/h}\,\intd x.
\end{equation}
Here, $\tilde{a}_h\in h^{0}\Cinfty(X\times B\times \la^*)$ is
the non-euclidean symbol of $\Op_{h}(a_h)$ defined in Lemma~\ref{Kap5-lemma-noneuclid-symb},
and $\psi$ is the phase function of \eqref{eq:psi-in-an-coords}.
\end{lemma}

\begin{proof}
We apply $\Op_h(a_h)$ to \eqref{eq:hi_h-off-diagonal}.
The rules for composing Schwartz kernels
justify interchanging the operator $\Op_h(a_h)$ with the integral (duality bracket).
In the notation of Lemma~\ref{Kap5-lemma-noneuclid-symb}, we get
\[
\Op_h(a_h)\phi_h(x)= \int_B \tilde{a}_h(x,b,\nu_h) e_{i\nu_h/h,b}(x)\, T_h(\intd b).
\]
Using the tensor product of distributions, we derive
\begin{align*}
\big(\Op_{h} & (a_h) \varphi_h\mid \chi\varphi'_h\big)_{L^2(X)} \\
 &= \int_X \int_{B^2} \chi(x) \tilde{a}_h(x,b,\nu_h) e_{i\nu_h/h,b}(x)
         e_{-w_0\cdot i\nu'_h/h,b'}(x) \, T_h(\intd b)\otimes \tilde T_h(\intd b') \,\intd x.
\end{align*}
We interchange the integral over $X$ with the duality bracket of distributions on $B^2$,
\[
\big(\Op_{h} (a_h) \varphi_h\mid \chi\varphi'_h\big)_{L^2(X)}
 = \int_{B^2} F_h(b,b',\nu_h,\nu'_h) \, T_h(\intd b)\otimes \tilde T_h(\intd b').
\]
Here we used $w_0^{-1}=w_0$,  and $-w_0\cdot\rho=\rho$.
\qed
\end{proof}

Consider the weight function 
$$d_h(gM,\nu,\nu') := d_{\nu/h,-w_o\cdot\nu'/h}(gM) =  e^{(i\nu/h+\rho)H(g)}  e^{(-i w_0\cdot\nu/h+\rho) H(gw_0)}.$$ 
Following \eqref{eq:Radon-trafo}, \eqref{eq:weight-function-off-diagonal}, and \eqref{repair}, we have the weighted Radon transform
\begin{gather*}
\mathcal R_h:\Ccinfty(G/M)\to \Ccinfty(G/MA)\subset\Cinfty(B^2),\\
\quad
(\mathcal R_hf)(gMA) =\int_A d_h(gaM,\nu_h,\nu'_h)f(gaM)\,\intd a,
\end{gather*}
and its dual $\mathcal R_h':\Dprime(B^2)\to\Dprime(G/M)$.
Further, \eqref{eq:PSGamma} suggests to define $PS^\Gamma_h:=PS^\Gamma_{\phi_h,\phi'_h}\in \Dprime(\Gamma\backslash G/M)$, so that
\[
\langle PS^\Gamma_h, f\rangle_{\Gamma\backslash G/M}
    = \langle PS^\Gamma_{\phi_h,\phi'_h}, f\rangle_{\Gamma\backslash G/M}
    =
    \langle {\mathcal R}_h' (T_h\otimes \tilde T_h), \chi f\rangle_{G/M}
\]
for $f\in\Ccinfty(\Gamma\backslash G/M)$, where $\chi\in\Ccinfty(G/M)$ is a smooth fundamental domain cutoff.

Given $\chi\in\Ccinfty(G/M)$ and $\chi_a\in\Ccinfty(\laos)$,
we define
\[ I_h=I_{h,\chi}:h^0\Ccinfty(T^*\XG)\to  h^{\dim N/2} C_c^\infty(G/M\times\laos)  \]
as follows.
For $S=\supp\chi$,
choose $\beta\in\Ccinfty(B^{(2)})\subset\Cinfty(B^2)$ as in Lemma~\ref{lem:beta}.
Denote by $\hat\beta$ the $A$-invariant lift of $\beta$ to $G/M$.
Recall the definition of the non-euclidean symbol $\tilde{a}_h\in h^{0}\Cinfty(G/M\times \la^*)$
of an operator $\Op_h(a_h)$ from Lemma~\ref{Kap5-lemma-noneuclid-symb}.
Following \eqref{def-Lh-off-diagonal} we set 
\begin{equation}
\label{def:op-Ih-off-diagonal}
\begin{aligned}
(I_h a_h) & (gM,\nu,\nu')  \\
   &:= \hat\beta(gM) \chi_a(\nu,\nu')
       \int_N \chi(gnM) \tilde{a}_h(gnM,\nu) e^{i(w_0\cdot\nu') H(n^{-1}w_0)/h} \intdbar n.
\end{aligned}
\end{equation}

We relate lifted quantum limits to Patterson--Sullivan distributions.

\begin{lemma}
\label{lem:qlim-and-PS-off-diagonal}
Set
\[
W_h(a)=\big(\Op_{\Gamma,h}(a)\varphi_h\mid \varphi'_h\big)_{L^2(\XG)},
\quad a\in\Ccinfty(T^* \XG).
\]
Assume that $\omega=\lim_h W_h$ in $\Dprime(T^* \XG)$ as $h\to 0$.
Assume further that $\lim_{h\to0}\nu_h=\theta$ and $\lim_{h\to0}\nu'_h=\theta'$ with $(\theta,\theta') \in\laos$.
Suppose $\chi$ is smooth fundamental domain cutoff,
and $\chi_a=1$ in a neighborhood of $(\theta,\theta')$.
Let $a_h\in S^{0,0}(T^* X)$ with principal symbol
$a=\lim_{h\downarrow 0} a_h\in\Ccinfty(T^*\XG)$.
Then, with $I_h=I_{h,\chi}$,
\begin{equation}
\label{mu-eq-limPS-off-diagonal}
\int_X a\intd \omega = \lim_{h\downarrow 0}
     \langle {\mathcal R}_h' (T_h\otimes \tilde T_h), (I_h a_h)(\cdot, \nu_h,\nu'_h)\rangle_{G/M}.
\end{equation}
\end{lemma}

\begin{proof}
Combine Proposition~\ref{prop:asymptotics-of-Fh-with-cutoff-off-diagonal}, Lemma~\ref{lem: mu probability measure}, and Lemma~\ref{lem:direct-calculation-off-diagonal}.
\qed
\end{proof}

\begin{remark}
\label{rem:technical-chi-ah-off-diagonal}
Observe that for any $\chi'\in\Ccinfty(G/M)$,
\[ \lim_{h\downarrow 0} \big(\Op_{h} (a_h) \varphi_h\mid \chi'\varphi'_h\big)_{L^2(X)} =0 \]
if $a_h\in S^{0,-1}(T^* X)$.
This observation will allow us to add terms to
\eqref{direct-calculation-off-diagonal} without changing the limit as $h\downarrow 0$.
\end{remark}

Recall from \eqref{eq:T*X-surjection} the $G$-equivariant map
$\Phi\colon G/M\times\la^*\to T^*X$.
If $\theta\in\la^*$ is regular, then $\Phi(\cdot,\theta):G/M\to T^* X$
is an imbedding having a joint level set as its range, \cite[Lemma 1.6]{Hi05}.
Since this map is proper, the push-forward of distributions,
\[
\Phi(\cdot,\theta)_*:\Dprime(\Gamma\backslash G/M)\to \Dprime(T^* \XG),
\]
is well-defined. Moreover, we can define an extension operator
\[
E_\theta:\Ccinfty(G/M)\to \Ccinfty(T^* X),
\quad
(E_\theta u)(\Phi(gM,\theta)) = u(gM).
\]

\begin{theorem}
\label{thm-W-PS-off-diagonal}
Let $(\varphi_h)_h,(\varphi'_h)_h\subset L^2(\XG)$ be sequences of normalized joint eigenfunctions, with purely imaginary spectral parameters $i\nu_h/h, i\nu'_h/h \in i\la^*$. Assume that $\omega=\lim_h W_h$ in $\Dprime(T^* \XG)$ as $h\to 0$. Assume further that $\lim_{h\to0}\nu_h=\theta$ and $\lim_{h\to0}\nu'_h=\theta'$ with  $(\theta,\theta') \in\laosplus$ such that
\begin{equation}
\label{assumption-on-lambda_h-off-diagonal}
\nu_h = \theta + O(h),\  \nu'_h = \theta' + O(h)\quad\text{as $h\downarrow 0$.}
\end{equation}
Then, with $\kappa$ defined in \eqref{eq:leading-msp-coeff},
\begin{equation}
\label{W-PS-off-diagonal}
\omega = \kappa(w_0\cdot\theta') \lim_{h\downarrow 0} (2\pi h)^{\dim N/2} \Phi(\cdot,\theta)_* PS^\Gamma_h
\quad\text{in $\Dprime(T^*\XG)$.}
\end{equation}
\end{theorem}

\begin{proof}
Let $a\in \Ccinfty(T^*\XG)$.
Let $a_h\in S^{0,0}_\Gamma$ with principal symbol $a=\lim_{h\downarrow 0} a_h$.
Applying Proposition~\ref{prop:asymptotics-of-Fh-with-cutoff-off-diagonal}, we obtain,
with $\chi$ now a smooth fundamental domain cutoff,
\begin{align*}
(I_h a_h)(gM,\nu,\nu')
  &= \kappa(w_0\cdot\nu') (2\pi h)^{\dim N/2} \big(\chi(gM) \tilde{a}_h(gM,\nu) + O(h)\big) \\
  &= \kappa(w_0\cdot\theta') (2\pi h)^{\dim N/2} \\
  &\quad \cdot\big(\chi(gM) \tilde{a}_h(gM,\theta) + O(|\nu-\theta|) +O(|\nu'-\theta'|)+ O(h)\big),
\end{align*}
in $h^{\dim N/2}\Ccinfty\big(G/M\times\laos\big)$.
Here we used Taylor expansion around $\theta$ for $\tilde{a}_h(gM,\nu)$ and Taylor expansion around $\theta'$ for $\kappa(w_0\cdot\nu')$.
Setting $\nu=\nu_h$ and $\nu'=\nu'_h$, and using the
assumption \eqref{assumption-on-lambda_h-off-diagonal}, we have, as $h\downarrow 0$,
\begin{align*}
(I_h a_h)(gM,\nu_h,\nu'_h)
  &= \kappa(w_0\cdot\theta') (2\pi h)^{\dim N/2} \big(\chi(gM) \tilde{a}_h(gM,\theta) + O(h)\big) \\
  &= \kappa(w_0\cdot\theta') (2\pi h)^{\dim N/2} \big(\chi(gM) a_h(\Phi(gM,\theta)) + O(h)\big) \\
  &= \kappa(w_0\cdot\theta') (2\pi h)^{\dim N/2} \big(\chi(gM) a(\Phi(gM,\theta)) + O(h)\big).
\end{align*}
The second equation follows from \eqref{noneuclid-symb-expansion}.

For $\ell>0$ sufficiently large, we shall modify $a_h$ by lower order terms, i.e.,
terms in $h^1\Ccinfty(T^*\XG)$, such that the above error term $O(h)$
gets replaced by $O(h^\ell)$.
This will, in view of Proposition~\ref{cor:Polynomial Bound}
and Lemma~\ref{lem:qlim-and-PS-off-diagonal}, imply
\[
\int_X a\intd\omega = \lim_{h\downarrow 0} \langle {\mathcal R}_h' (T_h\otimes \tilde T_h),\kappa(w_0\cdot\theta')(2\pi h)^{\dim N/2} \chi \Phi(\cdot,\theta)^* a\rangle,
\]
and hence the theorem.

Set $r_h(gM)= (I_h a_h)(gM,\nu_h,\nu'_h) - \kappa(w_0\cdot\theta') (2\pi h)^{\dim N/2} \chi(gM) a(\Phi(gM,\theta))$.
By the computation above,
\begin{equation}
\label{r_h-belongs-to-off-diagonal}
r_h\in h^{\ell+\dim N/2} \Ccinfty(G/M)
\end{equation}
with $\ell=1$.
Define
$a_h' = (2\pi h)^{-\dim N/2} E_\theta r_h \in h^\ell \Ccinfty(T^* X)$.
Choose $\chi'\in\Ccinfty(X)$ such that $\chi'=1$ on the support of $r_h$.
The computations above with $a_h$ replaced by $a_h'$ give
\begin{align*}
(I_{h,\chi'} a_h')(gM,\nu_h,\nu'_h)
  &= \kappa(w_0\cdot\theta') h^{\dim N/2} \big(\chi'(gM) a_h'(\Phi(gM,\theta)) + O(h^{1+\ell})\big) \\
  &= \kappa(w_0\cdot\theta') \big(r_h(gM) + O(h^{1+\ell})\big).
\end{align*}
We replace, in \eqref{mu-eq-limPS-off-diagonal}, $I_h a_h$ by $I_h a_h -\kappa^{-1} I_{h,\chi'} a_h'$.
By Remark~\ref{rem:technical-chi-ah-off-diagonal}, the formula \eqref{mu-eq-limPS-off-diagonal} remains true.
In addition, by the arguments above, we have a new remainder $r_h$ which
satisfies \eqref{r_h-belongs-to-off-diagonal} with $\ell$ replaced by $\ell+1$.
Arguing by induction over $\ell$, the proof follows.
\qed
\end{proof}

\begin{remark}
\begin{enumerate}
\item[(i)] If $\theta\neq\theta'$, then Remark~\ref{no-life-outside-the-diagonal}, combined with the method of non-stationary phase and Proposition~\ref{cor:Polynomial Bound} imply that $\omega=0$. Thus, in this case also the right hand side of \eqref{W-PS-off-diagonal} vanishes
\item[(ii)]  Combining Theorem~\ref{thm-W-PS-off-diagonal} for $\phi'_h=\phi_h$ with Remark~\ref{rem:A-equivariance-PS} yields yet another proof of the $A$-invariance of the lifted quantum limits (see Theorem~\ref{Kap5-Prop-supp-invar-Wigner}, where $|\nu_h|=1$).

\end{enumerate}
\end{remark}

%% file: hhs-arXiv.bbl
\begin{thebibliography}{}



\bibitem{AS} N. Anantharaman, L. Silberman, \textit{A Haar Component for Quantum Limits in Locally Symmetric Spaces}, preprint, 2010, arXiv:1009.492.

\bibitem{AZ} N. Anantharaman, S. Zelditch, \textit{Patterson--Sullivan Distributions and Quantum Ergodicity}, Ann. Henri Poincar\'{e} {\bf 8} (2007), 361--426.

\bibitem{AZ2} N. Anantharaman, S. Zelditch, \textit{Intertwining the geodesic flow and the Schr\"odinger group on hyperbolic surfaces}, preprint, 2010, arXiv:1010.0867.

\bibitem{An91} J.-P. Anker, \textit{A Basic Inequality for Scattering Theory of Riemannian Symmetric Spaces of the Noncompact Type}. Amer.~J.~Math.. {\bf 113} (1991), 391--398.

\bibitem{BS87} E. van den Ban, H. Schlichtkrull, \textit{Asymptotic expansions and boundary values of eigenfunctions on Riemannian symmetric spaces}, J. reine angew. Math. {\bf 380} (1987), 108--165.



\bibitem{DS} M.~Dimassi and J.~Sj{\"o}strand, \textit{Spectral Asymptotics in the Semi-Classical Limit}, Cambridge University Press, Cambridge (1999).

\bibitem{DKV1} J. J. Duistermaat, J. A. C. Kolk, V. S. Varadarajan, \textit{Spectra of compact locally symmetric manifolds of negative curvature}, Invent. Math. {\bf 52} (1979), 27--93.

\bibitem{DKV} J. J. Duistermaat, J. A. C. Kolk, V. S. Varadarajan, \textit{Functions, Flows and oscillatory integrals on flag manifolds and conjugacy classes in real semisimple Lie groups}, Compositio Math. {\bf 49} (1983), 309--398.


\bibitem{EZ03} L.E. Evans, M. Zworski \textit{Lectures on Semiclassical Analysis}, UC Berkeley. 

\bibitem{Ha10} S.~Hansen, \textit{Rayleigh-type surface quasimodes in general linear elasticity}, preprint 2010, arXiv:1008.2930.


\bibitem{Heckman} G.~Heckman, \emph{Projection of Orbits and Asymptotic Behaviour of Multiplicities for Compact Lie Groups}, Dissertation, Univ. Utrecht, 1980.



\bibitem{He94} S. Helgason, \textit{Geometric Analysis on Symmetric Spaces}, Mathematical surveys and monographs, American Mathematical Society, Providence, RI (1994).

\bibitem{He00} S. Helgason, \textit{Groups and Geometric Analysis}, Mathematical surveys and monographs, American Mathematical Society, Providence, RI (2000).

\bibitem{He01} S. Helgason, \textit{Differential Geometry, Lie Groups, and Symmetric Spaces}, Graduate Studies in Mathematics, American Mathematical Society, Providence, RI (2001).

\bibitem{HS09} J. Hilgert, M. Schr\"oder, \textit{Patterson--Sullivan Distributions for Rank One Symmetric Spcaces of the Non-Compact Type}. arXiv0909.2142

\bibitem{Hi05} J. Hilgert, \textit{An Ergodic Arnold--Liouville Theorem for Locally Symmetric Spaces}, ``Twenty Years of Bialowieza: A Mathematical Anthology''. S.T. Ali et al. eds., World Scientific, Singapore (2005).

\bibitem{Horm} L.~H{\"o}rmander, \textit{The Analysis of Linear Partial Differential Operators}, volumes I-IV, Springer-Verlag, Berlin and New York (1983--1985).


\bibitem{KKMOOT} M.~Kashiwara, A.~Kowata, K.~Minemura, K.~Okamoto, T.~{\=O}shima, and M.~Tanaka. \textit{Eigenfunctions of invariant differential operators on a symmetric space}, Ann.~of Math. {\bf 107} (1978), 1--39.



\bibitem{Schlichtkrull} H. Schlichtkrull, \textit{Hyperfunctions and Harmonic Analysis on Symmetric Spaces}, Progress in Math. {\bf 49}, Birkh\"auser, 1984.

\bibitem{S} M.~Schr{\"o}der, \emph{Patterson-Sullivan distributions for rank one symmetric spaces of the noncompact type}, Dissertation, Univ. Paderborn, 2010.

\bibitem{Sh05a} V.A.~ Sharafutdinov, \emph{Geometric symbol calculus for pseudo differential operators I and II}. Siber. Adv. Math.~{\bf 15}, no.3~(2005), 81--125, and no.4~(2005), 71--95.


\bibitem{SV07} L.~Silberman, A. Venkatesh, \emph{On quantum unique ergodicity for locally symmetric spaces}, Geom.~Funct.~Anal.~{\bf 17} (2007), 960--998.

\bibitem{T} F. Treves, \textit{Topological Vector Spaces, Distributions and Kernels}, Acad. Press, (1967)

\bibitem{Tay81} M. E. Taylor, \textit{Pseudodifferential Operators}, Princeton Univerity Press, Princeton, New Jersey, (1981).



\bibitem{Wal2} N. R. Wallach, \textit{Real reductive groups 1}, Academic Press, Pure and Applied Mathematics, San Diego (1988).

\bibitem{Wi} H. Widom, \textit{A complete symbolic calculus for pseudodifferential operators.} Bull. Sci. Math (2) {\bf 104} (1980), 19--63.

\bibitem{Wil} F. L. Williams, \textit{Lectures on the spectrum of $L^2(\Gamma\backslash G)$}, Pitman Research Notes in Mathematics Series 242, Essex (1991).

\bibitem{Wol} S. Wolpert, \textit{The Modulus of Continuity for $\Gamma_0(m)$ Semi-Classical Limits}, Comm. Math. Phys. {\bf 216} (2001), 313--323.

\bibitem{Z84} S. Zelditch, \textit{Pseudo-differential Analysis on Hyperbolic Surfaces}, J. of Funct. Anal. {\bf 68} (1986), 72--105.

\bibitem{Z87} S. Zelditch, \emph{Uniform distribution of eigenfunctions on compact hyperbolic surfaces}, Duke Math. J. {\bf 55} (1987), 919--941.

\bibitem{Z89} S. Zelditch, \emph{The averaging method and ergodic theory for pseudo-differential operators on compact hyperbolic surfaces}, J. of Funct. Anal. {\bf 82} (1989), 38--68.

\bibitem{Zel09} S. Zelditch, \textit{Local and global analysis of eigenfunctions on Riemannian manifolds}, preprint (2009), arXiv:0903.3420v1

\end{thebibliography}
